\documentclass[]{article}
\usepackage{graphicx}
\usepackage{graphics}
\usepackage{amssymb,amsmath,amsthm,epsfig}
\newcommand{\J}{P^{(\alpha,\beta)}}
\newcommand{\wJ}{\widetilde P^{(\alpha,\beta)}}
\newtheorem{theorem}{Theorem}

\newtheorem{corollary}{Corollary}
\newtheorem{definition}{Definition}
\newtheorem{proposition}{Proposition}

\newtheorem{remark}{Remark}
\newcommand{\norm}[1]{{\left\|{#1}\right\|}}
\newcommand{\R}{\mathbb R}

\newcounter{reh}
\newcounter{rek}
\begin{document}

\begin{center}
{\large {\bf Generalized Prolate Spheroidal Wave Functions: Spectral Analysis  and Approximation 
of Almost Band-limited Functions.}}\\
\vskip 1cm Abderrazek Karoui$^a$  {\footnote{
Corresponding author: Abderrazek Karoui, Email: abderrazek.karoui@fsb.rnu.tn
This work was supported by the DGRST research Grant UR13ES47.}} and Ahmed Souabni$^a$
\end{center}
\vskip 0.5cm {\small

\noindent $^a$ University of Carthage,
Department of Mathematics, Faculty of Sciences of Bizerte, Tunisia.
}\\

\noindent{\bf Abstract}--- In this work, we first give various explicit and local estimates
of the eigenfunctions of a  perturbed  Jacobi differential operator. These eigenfunctions
generalize the famous classical  prolate spheroidal wave functions (PSWFs), founded in 1960's by D. Slepian and his co-authors and corresponding to the case $\alpha=\beta=0.$ They also generalize the new  PSWFs
introduced and studied  recently in \cite{Wang2}, denoted by GPSWFs and corresponding to the case $\alpha=\beta.$ 
The main content of this work is devoted  to the previous interesting special case $\alpha=\beta >- 1.$ In particular, we give further computational improvements, as well as some  useful explicit and local  estimates of the GPSWFs. More importantly, by using the concept of a restricted Paley-Wiener space,  we relate the GPSWFs to the  solutions of a generalized energy maximisation problem. As a consequence, many
desirable spectral properties of the self-adjoint compact integral operator associated with the GPSWFs are deduced 
 from the rich literature  of the  PSWFs. In particular,  we show that 
the GPSWFs are well adapted for the spectral approximation of the classical $c-$band-limited as well as almost $c-$band-limited functions.  Finally, we provide the reader with some numerical examples that illustrate the different results of this work.\\

\noindent {2010 Mathematics Subject Classification.} Primary  42C10, 33E10,  Secondary 34L10, 41A30..\\

\noindent {\it  Key words and phrases.}  Sturm-Liouville operators, finite weighted  Fourier transform, eigenvalues and eigenfunctions, special functions, prolate spheroidal wave functions,
band-limited functions. \\

\section{Introduction}

We first recall that for a bandwidth $c>0,$ the infinite countable set of the eigenfunctions of the finite Fourier transform
${\displaystyle \mathcal F_c,}$ defined on $L^2(-1,1)$ by ${\displaystyle \mathcal F_c f(x)=\int_{-1}^1 e^{icxy} f(y)\, dy,}$ are known as the prolate spheroidal wave functions (PSWFs). They have been extensively studied in the literature,
since the pioneer work on the subject  of D. Slepian and his collaborators, see \cite{Landau, Slepian1, Slepian2}.
The interest in the PSWFs  is essentially due  to their  wide range of  applications in different 
scientific research area such as signal processing, physics, applied mathematics, see for example \cite{Boyd1, Boyd2,
Hogan, Wang}.  Recently, in \cite{Wang2}, the authors have 
given a generalization of the PSWFs by considering the eigenfunctions of the special case of the weighted Fourier transform 
$\mathcal F_c^{(\alpha)},$ defined by ${\displaystyle \mathcal F_c^{(\alpha)} f(x)=\int_{-1}^1 e^{icxy}  f(y)\,(1-y^2)^{\alpha}\, dy,\, \alpha > -1.}$ Note that  although, the extra weight function $\omega_{\alpha}(x) = (1-x^2)^{\alpha},$ generates new computational complications,
the resulting eigenfunctions have some advantages over the classical PSWFs. In this paper, we first give some useful analytic and local estimates
of a more general Jacobi-type PSWFs. This is done by using special spectral techniques from the  theory of Sturm-Liouville operators applied
to the following Jacobi perturbed differential operator, defined on $C^2([-1,1])$ by 
\begin{equation}\label{Sturm_operator}
\mathcal L^{(\alpha,\beta)}_c \varphi(x)= (1-x^2)\varphi''(x)+\left((\beta-\alpha)-(\alpha+\beta+2)x\right)\varphi'(x) - c^2 x^2 \varphi(x)=
\mathcal L_0 \varphi(x)-c^2 x^2 \varphi(x).
\end{equation}
Note that in the limiting case $c=0,$ the eigenfunctions of the previous differential operator are reduced to the  Jacobi polynomials
$P_k^{(\alpha,\beta)}.$ Moreover, in the special case $c>0,$ $\alpha=\beta =0,$ these eigenfunctions correspond to  the classical Slepian prolate spheroidal wave functions. Moreover, in the case where $\alpha=\beta>-1,$
it has been shown in \cite{Wang2} that the finite weighted Fourier transform $\mathcal F_c^{(\alpha)}$ commutes with 
$\mathcal L^{(\alpha,\alpha)}_c.$ Hence,  both operators have the same eigenfunctions,
called generalized prolate spheroidal wave functions (GPSWFs) and  simply  denoted by $\psi^{(\alpha)}_{n,c},\, n\geq 0.$ They are solutions  of the following integral equation

\begin{equation} \label{eq1}
\mathcal F_c^{(\alpha)} \psi^{(\alpha)}_{n,c} (x) = 
\int_{-1}^{1} e^{icxy} \psi^{(\alpha)}_{n,c}(y) \omega_{\alpha}(y)\, dy =\mu_{n}^{(\alpha)}(c) \psi^{(\alpha)}_{n,c}(x),\quad |x|\leq 1.
\end{equation}
Here, $\mu_{n}^{(\alpha)}(c)$ is the eigenvalue of the integral operator $\mathcal F_c^{(\alpha)},$ associated with the 
eigenfunction $ \psi^{(\alpha)}_{n,c}.$

The  important part of this work is devoted to the study of the GPSWFs, that is  $\alpha=\beta.$ 
In particular, we give their  analytic extensions  to the whole real line. As a consequence,
we obtain an explicit and practical  formula (in terms of a ratio of two fast converging series) for computing the eigenvalues $\mu_n^{(\alpha)}(c).$ Note that the behaviour as well
as the decay  rate of these eigenvalues, play a crucial role in most  applications of the GPSWFs.
In \cite{Wang2}, by using an heuristic asymptotic
analysis, the authors have given an asymptotic super-exponential decay rate of the  $|\mu_n^{(\alpha)}(c)|.$  In the second part of this work,
we prove that for any $\alpha \geq 0,$ the  super-exponential decay rate of $|\mu_n^{(\alpha)}(c)|$ starts holding from the  plunge region around 
${\displaystyle n = \frac{2c}{\pi}.}$  The proof of  this result is based on the  characterization of the  GPSWFs as solutions of a generalized energy maximization problem, over a restricted Paley-Wiener space $B_c^{(\alpha)},$ given by 
$$ B_c^{(\alpha)}=\{ f\in L^2(\mathbb R),\,\, \mbox{Support } \widehat f\subset [-c,c],\,\, \widehat f \in 
L^2((-c,c),\omega_{-\alpha}(\frac{\cdot}{c})).   \}$$ 
Here, $\widehat f$ denotes the Fourier transform of $f\in L^2(\mathbb R),$ defined by
$$\widehat f(\xi)= \lim_{A\rightarrow +\infty}\int_{-A}^A e^{-ix \xi} f(x)\, dx,\quad \xi \in \mathbb R .$$

More precisely, for a real number $\alpha >-1,$ let  $J_{\alpha}$ denote the Bessel function of the first type and order 
$\alpha$ and consider the self-adjoint compact operator ${\displaystyle \mathcal Q_c^{\alpha}=\frac{c}{2\pi}
\mathcal F_c^{{\alpha}^*} \circ \mathcal F_c^{\alpha},}$ defined on $L^2{(I, \omega_{\alpha})},\, I=[-1,1]$ by 
\begin{equation}
\mathcal Q_c^{\alpha} g (x) = \int_{-1}^1 \frac{c}{2 \pi}\mathcal K_{\alpha}(c(x-y)) g(y) \omega_{\alpha}(y) \, dy,
\quad \mathcal K_{\alpha}(x)=\sqrt{\pi} 2^{\alpha+1/2}\Gamma(\alpha+1) \frac{J_{\alpha+1/2}(x)}{x^{\alpha+1/2}}.
\end{equation}
By  rewriting the energy  maximization problem in term of the previous integral operator,
one gets a characterization of the eigenvalues ${\displaystyle \lambda_n^{(\alpha)}(c)= \frac{c}{2\pi } |\mu_n^{(\alpha)}(c)|^2}$ of $\mathcal Q_c^{\alpha}$ as a countable sequence generated by the energy problem. From this, we conclude that the $\lambda_n^{(\alpha)}(c)$ decay with respect 
to the parameter $\alpha,$ that is for $c>0,$ and any $n\in \mathbb N,$ 
$$ 0 <  \lambda_n^{(\alpha)}(c) \leq  \lambda_n^{(\alpha')}(c) < 1,\quad \forall \, \alpha \geq \alpha' \geq 0.$$
 Hence, by using the precise behaviour as well as the sharp decay rate of the $\lambda_n^{(0)}(c),$ given in 
\cite{Bonami-Karoui3}, one gets a similar behaviour and decay rate  for the $(\lambda_n^{(\alpha)}(c))_n,$ for any $\alpha \geq 0.$\\

This work is organized as follows. In section 2, we   give some mathematical preliminaries on Jacobi polynomials and their finite weighted Fourier transform. Also, we describe the Bouwkamp method for the Jacobi series expansion of the eigenfunctions $\psi_{n,c}^{(\alpha,\beta)}$ of the differential operator $\mathcal L_c^{(\alpha,\beta)}.$ Then, we give some explicit and local estimates of these eigenfunctions. These estimates will be particularly useful in the subsequent study of the GPSWFs, given by sections 3 and 4. In section 3, we first give an improvement as well as a new kind of decay rate of the  Gegenbauer's series expansion coefficients of the GPSWFs $\psi_{n,c}^{(\alpha)}.$ Then, we give the analytic extension of the GPSWFs to the whole real line, as well as an explicit and practical  formula ( as a ratio of two fast converging series) for the accurate computation of the eigenvalues $\mu_n^{(\alpha)}(c)$ and consequently of $\lambda_n^{(\alpha)}(c).$ In section 4, we characterize the GPSWFs as solutions of a generalized
energy maximization problem and we prove the monotonicity of the $\lambda_n^{(\alpha)}(c)$ with respect to the parameter $\alpha.$ Moreover, we show that the GPSWFs are well adapted for the approximation of the classical $c-$band-limited as well as almost $c-$band-limited functions. Finally, in section 5, we provide the reader with numerical examples that illustrate the different results of this work.\\

\noindent
{\bf Notations and normalizations:} The following notations will be frequently used in this work, 
$$\omega_{\alpha,\beta}(y)=(1-y)^{\alpha}(1+y)^{\beta},\,\, \quad \omega_{\alpha}(y)= (1-y^2)^{\alpha}, \,\, 
I=[-1,1].$$
Moreover, the eigenfunctions $\psi_{n,c}^{(\alpha,\beta)}$ and the GPSWFs $\psi_{n,c}^{(\alpha)}$ are normalized 
so that 
$$\int_{-1}^1 \left(\psi_{n,c}^{(\alpha,\beta)}(t)\right)^2 \, \omega_{\alpha,\beta}(t)\, dt =1,\quad 
\int_{-1}^1 \left(\psi_{n,c}^{(\alpha)}(t)\right)^2 \, \omega_{\alpha}(t)\, dt =1.$$

\section{Eigenfunctions of a perturbed Jacobi differential operator.} 

In this section, we  give a description  of the series expansion of the eigenfunctions $\psi_{n,c}^{(\alpha,\beta)}(x)$
of $\mathcal L_c^{(\alpha,\beta)},$  given by \eqref{Sturm_operator} and with respect to  the basis of normalized 
Jacobi polynomials $\mathcal B=\{\wJ_k,\, k\geq 0\}.$  Also, we give some properties as well as local estimates 
of $\psi_{n,c}^{(\alpha,\beta)}(x),$ generalizing some of those given in \cite{Bonami-Karoui1, Bonami-Karoui2} in the special case $\alpha=\beta=0.$ For this purpose, we  need the following mathematical preliminaries.
 
\subsection{Mathematical preliminaries} 

We first recall that for two real numbers $\alpha, \beta >-1,$ the Jacobi polynomials $\J_k$ are given by the following three
term recursion formula
\begin{equation}\label{recursion1}
\J_{k+1}(x)= (A_k x + B_k)\J_{k}(x) -C_k \J_{k-1}(x),\quad x\in [-1,1].
\end{equation}
with  $\J_0(x)=1,\quad \J_1(x)=\frac{1}{2}(\alpha+\beta+2)x +\frac{1}{2}(\alpha+\beta).$  Here, 
\begin{eqnarray}\label{recursion2}
A_k &=&\frac{(2k+\alpha+\beta+1)(2k+\alpha+\beta+2)}{2 (k+1)(k+\alpha+\beta+1)},\quad B_k=\frac{(\alpha^2-\beta^2)(2k+\alpha+\beta+1)}{2(k+1)(k+\alpha+\beta+1)(2k+\alpha+\beta)}\nonumber\\
C_k&=&\frac{(k+\alpha)(k+\beta)(2k+\alpha+\beta+2)}{(k+1)(k+\alpha+\beta+1)(2k+\alpha+\beta)}.
\end{eqnarray}
In the sequel, we let $\wJ_k$ denote the normalized Jacobi polynomial of degree $k$ so that 
$$ \| \wJ_k \|^2_{L^2_{\omega_{\alpha,\beta}(I)}}=\int_{-1}^1 (\wJ_k(y))^2 \omega_{\alpha,\beta}(y)\, dy =1.$$
 It is well known that in this case, we have 
\begin{equation}\label{JacobiP}
\wJ_{k}(x)= \frac{1}{\sqrt{h_k}}\J_k(x),\quad h_k=\frac{2^{\alpha+\beta+1}\Gamma(k+\alpha+1)\Gamma(k+\beta+1)}{k!(2k+\alpha+\beta+1)\Gamma(k+\alpha+\beta+1)}.
\end{equation}
Straightforward computations give us the following useful identities
\begin{equation}\label{recursion3}
\wJ_{k+1}(x)= (a_k x + b_k)\wJ_{k}(x) -c_k \wJ_{k-1}(x),
\end{equation}
\begin{equation}\label{coefficients}
a_k=\sqrt{\frac{h_{k}}{h_{k+1}}} A_k,\quad b_k=\sqrt{\frac{h_{k}}{h_{k+1}}} B_k,\quad  c_k=\sqrt{\frac{h_{k-1}}{h_{k+1}}} C_k.
\end{equation}
\begin{eqnarray}\label{Eq2.1}
x^2 \wJ_k(x) &=& \frac{1}{a_k a_{k+1}}\wJ_{k+2}(x)-\Big(\frac{b_{k+1}}{a_k  a_{k+1}}+\frac{b_k}{a_k^2}\Big)\wJ_{k+1}(x)
+\Big(\frac{c_{k+1}}{a_k a_{k+1}}+\frac{b_k^2}{a_k^2} 
+\frac{c_k}{a_k a_{k-1}}\Big)\wJ_k(x)\nonumber \\
&&-\Big(\frac{c_k b_k}{a_k^2}+\frac{c_k b_{k-1}}{a_k a_{k-1}}\Big)\wJ_{k-1}(x) + \frac{c_k c_{k-1}}{a_k a_{k-1}}\wJ_{k-2}(x)
\end{eqnarray}
The explicit expressions and  bounds of the different moments of the weight function $\omega_{\alpha,\beta}$
as well as of  the Jacobi polynomials $\wJ_k$ will be frequently needed in this work.  For this purpose, 
we first recall the following useful inequalities for the Gamma function, see \cite{Batir}, 
\begin{equation}\label{Ineq2}
\sqrt{2e} \left(\frac{x+1/2}{e}\right)^{x+1/2}\leq \Gamma(x+1)\leq \sqrt{2\pi} \left(\frac{x+1/2}{e}\right)^{x+1/2},\quad x>0. 
\end{equation}
Next, for an integer $k\geq 0,$  let 
\begin{equation}
\label{moment_weight}
I_k^{\alpha,\beta}=\int_{-1}^1 y^{k}\omega_{\alpha,\beta}(y)\, dy
\end{equation}
be the $k-$th moment of $\omega_{\alpha,\beta}.$ To get an upper bound for  $I_k^{\alpha,\beta},$
we may assume that $\alpha\geq \beta.$ In this case, we have 
\begin{eqnarray}\label{Ineq4}
I_k^{\alpha,\beta} &=& \int_0^1 y^{k} (1-y)^{\alpha}(1+y)^{\beta}\, dy + \int_0^1 y^{k} (1-y)^{\beta}(1+y)^{\alpha}\, dy 
\leq  \int_0^1 y^{k} (1-y^2)^{\alpha}\, dy \nonumber\\
&& + 2^{\alpha-\beta} \int_0^1 y^{k} (1-y^2)^{\beta}\, dy
\leq  \frac{1}{2} \left( B(k/2+1/2,\alpha+1)+2^{\alpha-\beta} B(k/2+1/2,\beta+1)\right).
\end{eqnarray}
Here $B(\cdot,\cdot)$ is the Beta function given by ${\displaystyle B(x,y)=\frac{\Gamma(x) \Gamma(y)}{\Gamma(x+y)},\, x, y >0.}$
Moreover, by using \eqref{Ineq2}, taking into account that for any real number $a>-1,$ the function 
\begin{equation}
\label{function}
\varphi(x)=  \left(1+\frac{a}{x}\right)^{a+x}, \quad x\geq 1
\end{equation}
 is decreasing  on $[1,\infty)$ to $e^a$ and by using some straightforward computations, one gets 
$$I_k^{\alpha,\beta}\leq \frac{1}{2}\sqrt{\frac{\pi}{e}}\frac{2^{1+\beta}}{k^{1+\beta}} \left(\Gamma(1+\alpha)+2^{\alpha-\beta }\Gamma(1+\beta)\right),\quad \alpha\geq \beta.$$
Consequently, for any real numbers $\alpha,\beta >-1,$ we have 
\begin{equation}\label{Ineq5}
I_k^{\alpha,\beta} 
\leq \frac{C_{\alpha,\beta}}{k^{1+\min(\alpha,\beta)}},\quad\forall\, \alpha, 
\beta > -1,\qquad C_{\alpha,\beta}=\left(2^{\alpha}+2^{\beta }\right)\sqrt{\frac{\pi}{e}}\Gamma(1+\max(\alpha,\beta)).
\end{equation}
In the special case where $\alpha =\beta,$ and by using the parity of $\omega_{\alpha}(y)=(1-y^2)^{\alpha}$ as well as the previous  bound, one gets
\begin{equation}\label{Ineqq5}
I_{2k+1}^{\alpha,\alpha}=0,\quad I_{2k}^{\alpha,\alpha} = B(k+1/2,\alpha+1)\leq 
\sqrt{\frac{\pi}{e}}\frac{\Gamma(\alpha+1)}{k^{1+\alpha}}.
\end{equation}
Also, note that for given integers $k\geq n\geq 0,$ and by using the 
Rodrigues formula for the Jacobi polynomials, one gets the following formula for the $k-$th moments of $\wJ_n,$
with $k\geq n,$  
\begin{equation}\label{moment_Jacobi}
M_{k,n}=\int_{-1}^1 x^k \wJ_n(x)\omega_{\alpha,\beta}(x)\, dx 
=\frac{1}{2^n\sqrt{h_n}}{k\choose n}\int_{-1}^1 x^{k-n}(1-x)^{n+\alpha}(1+x)^{n+\beta}\, dx.
\end{equation}
In particular, if $\alpha=\beta,$ one gets  
\begin{equation}\label{moment_Jacobi2}
M_{k,n}=\int_{-1}^1 x^k \widetilde P_n^{(\alpha,\alpha)}(x)\omega_{\alpha}(x)\, dx 
=\left\{ \begin{array}{ll} 0 & \mbox{ if } k< n \mbox{ or } k-n \mbox{ is odd}\\
\frac{{k\choose n}}{2^n\sqrt{h_n}} B\left(\frac{k-n+1}{2},n+\alpha+1\right)& \mbox{ otherwise.}\end{array}\right.
\end{equation}

On the other hand, it is  interesting to note that the weighted finite Fourier transform of Jacobi polynomial is given by the following explicit expression, see [\cite{NIST}, p.456],
\begin{equation}\label{fourier_Jacobi}
\int_{-1}^{1} e^{ixy} P_k^{(\alpha , \beta)}(y) \omega_{\alpha,\beta}(y)\,  dy=
\frac{(ix)^k e^{ix}}{k!}2^{k+\alpha+\beta+1}B(k+\alpha+1,k+\beta+1) \\
_1F_1(k+\alpha+1,2k+\alpha+\beta+2,-2ix),
\end{equation}
where $B(x,y)$ is the Beta function and $_1 F_1(a,b,c)$ is the Kummer's function. It is well known,
see [\cite{NIST}, p.326] that the Kummer's function
has the following integral representation
\begin{equation}\label{Kummer1}
_1F_1(a,b;z)=\frac{1}{B(a,b)}\int_0^1 e^{zt}t^{a-1}(1-t)^{b-a-1}dt,\quad z\in \mathbb C,\quad \mathcal Re(b)> \mathcal Re(a)>0.
\end{equation}

\subsection{Computation and first properties of  the eigenfunctions of $\mathcal L_c^{(\alpha,\beta)}.$}

In this paragraph, we first  describe the  Bouwkamp method for the computation of the bounded eigenfunctions and 
the corresponding eigenvalues of the operator  $\mathcal L_c^{(\alpha,\beta)},$ given by \eqref{Sturm_operator}. Then, we give some general properties of 
these eigenfunctions. Note that Bouwkamp method can be briefly described as  the representation of a perturbed version of classical orthogonal polynomials
differential operator. This representation is done by the use of the original  classical orthogonal polynomials. 
In our case, we consider  the Jacobi  orthonormal basis of $L^2(I,\omega_{\alpha,\beta}),$ given by $\mathcal B^{\alpha,\beta}=\{ \wJ_k(x),\, k\geq 0\}.$ Then,  thanks to this method,
the computation of the bounded eigenfunctions $\psi^{(\alpha,\beta)}_{n,c}$ of $\mathcal L_c^{(\alpha,\beta)}$ and their associated
eigenvalues $\chi_n(c)$ is reduced to the  computation of the eigenvectors and the  associated eigenvalues of the infinite order
matrix representation of  $\mathcal L_c^{(\alpha,\beta)}$ with respect to the basis $\mathcal B^{\alpha,\beta}.$
It is interesting to note that only a finite number of the main diagonals of this representation matrix  are not identically zeros. To the best of our knowledge,  C. Niven, was the first to use this method in the early 1880's, see \cite{Niven}.

Note that  since  $\psi^{(\alpha,\beta)}_{n,c}\in L^2(I, \omega_{\alpha,\beta}),$ then its series expansion
with respect to the basis $\mathcal B^{\alpha,\beta}$ is given by 
\begin{equation}\label{expansion1}
\psi^{(\alpha,\beta)}_{n,c}(x) =\sum_{k\geq 0} \beta_k^n \wJ_k (x),\quad x\in [-1,1].
\end{equation}
By combining \eqref{expansion1} and the facts that   
$$-\mathcal L_c^{(\alpha,\beta)} \psi^{(\alpha,\beta)}_{n,c}(x)=\chi_n(c) \psi^{(\alpha,\beta)}_{n,c}(x),\quad -\mathcal L_0^{(\alpha,\beta)} \wJ_k(x)=
\chi_k(0)\wJ_k(x),\quad \chi_k(0)=k(k+\alpha+\beta+1),$$ one can easily check that the expansion coefficients $(\beta_k^n)_{k\geq 0},\, n\geq 0$ and the 
eigenvalues $(\chi_n(c))_{n\geq 0}$ are given by the following infinite order eigensystem
\begin{equation}\label{eigensystem1}
\mathbf D^{\alpha,\beta} \cdot \mathbf B_n = \chi_n(c)\,  \mathbf B_n,\quad \mathbf D^{\alpha,\beta} = \left[ d_{i,j}\right]_{i,j\geq 0},\quad 
\mathbf B_n = [\beta_k^n, k\geq 0]^T.
\end{equation}
Here, $\mathbf D^{\alpha,\beta}$ is a $5-$diagonals  matrix representation of the operator $-\mathcal L_c^{(\alpha,\beta)}$ with coefficients 
given by  
\begin{eqnarray}\label{eigensystem2}
d_{i,i-2}&=& d_{i,i-2} = c^2\frac{1}{a_{i-1}a_{i-2}},\quad d_{i,i-1}=-c^2\left(\frac{b_i}{a_ia_{i-1}}+\frac{b_{i-1}}{a_{i-1}^2}\right)\nonumber \\
d_{i,i}&=&i(i+\alpha+\beta+1)+c^2\left(\frac{c_{i+1}}{a_i a_{i+1}}-\frac{b_i^2}{a_i^2}+\frac{c_{i}}{a_i a_{i-1}}\right),\,
d_{i,i+1}= c^2\left( \frac{c_{i+1} b_{i+1}}{a_{i+1}^2}+\frac{c_{i+1} b_i}{a_{i+1} a_i}\right),\nonumber\\
 d_{i,j}&=& 0,\,\,\,\mbox{ if } |j-i|\geq 3.
\end{eqnarray}
We recall that the coefficients $a_i, b_i, c_i$ are given by \eqref{recursion2} and \eqref{coefficients}. In the special case where 
$\alpha =\beta,$ we have $b_i=0,$ so that the previous eigensystem is reduced to a symmetric tri-diagonal system. In this case, for a fixed
integer $n\geq 0,$ the sequence $(\beta_k^n)_{k\geq 0}$ satisfies the following eigensystem
\begin{equation}\label{eigensystem3}
c^2\frac{1}{a_{k+2}a_{k+1}} \beta_{k+2}^n + \left(k (k+2\alpha+1)+c^2 \left(\frac{c_{k+1}}{a_k a_{k+1}}+\frac{c_{k}}{a_k a_{k-1}}\right)\right)
\beta_k^n
+c^2\frac{1}{a_{k}a_{k-1}} \beta_{k-2}^n= \chi_n(c) \beta_k^n.
\end{equation}
An expanded form of this system is given by
\begin{eqnarray}\label{eigensystem}
\lefteqn{\frac{\sqrt{(k+1)(k+2)(k+2\alpha+1)(k+2\alpha+2)}}{(2k+2\alpha+3)\sqrt{(2k+2\alpha+5)(2k+2\alpha +1)}} c^2 \beta_{k+2}^n
 + \big( k(k+2\alpha+1)+c^2 \frac{2k(k+2\alpha+1)+2\alpha-1}{(2k+2\alpha+3)(2k+2\alpha-1)} \big)
\beta_k^n}\nonumber  \\
&&\hspace*{2cm} + \frac{\sqrt{k(k-1)(k+2\alpha)(k+2\alpha-1)}}{(2k+2\alpha-1)\sqrt{(2k+2\alpha+1)(2k+2\alpha-3)}} c^2
\beta_{k-2}^n= \chi_n(c) \beta_k^n, \quad k\geq 0.
\end{eqnarray}

The following proposition provides us with some  properties of the eigenfunctions $\psi^{(\alpha,\beta)}_{n,c}(x)$ and eigenvalues
$\chi_n(c),$  generalizing some known properties for Jacobi polynomials.

\begin{proposition}\label{properties1}
For given real numbers $c>0,$ $\alpha, \beta > -1,$ let $\psi_{n,c}^{(\alpha,\beta)}$  be the   $n-$th eigenfunction
associated with  $\mathcal L_c^{(\alpha,\beta)},$ and  
normalized so that $\|\psi_{n,c}^{(\alpha,\beta)}\|_{L^2(I,\omega_{\alpha,\beta})}=1.$
Then we have \\
$(P_1)$ The set $\mathcal B=\{ \psi_{n,c}^{(\alpha,\beta)},\, n\geq 0 \}$ is an
orthonormal basis of $L^2(I,\omega_{\alpha,\beta}).$\\
$(P_2)$ If $\psi_{n,c}^{(\beta,\alpha)}$ is the $n-$th normalized eigenfunction of $\mathcal L_c^{(\beta,\alpha)},$ then
 $\psi_{n,c}^{(\alpha,\beta)}$ and $\psi_{n,c}^{(\beta,\alpha)}$ are associated to the same eigenvalue $\chi_n(c).$ Moreover, they 
 are related  to each others by the 
following rule
\begin{equation}
\label{symmetry}
\psi_{n,c}^{(\alpha,\beta)}(-x)= (-1)^n \psi_{n,c}^{(\beta,\alpha)}(x),\quad x\in \mathbb R.
\end{equation}
\end{proposition}

\noindent
{\bf Proof:}  Property $(P_1)$ follows from the general spectral theory of Sturm-Liouville operators. To prove $(P_2),$
we use the following  well known property for Jacobi polynomials, see [\cite{Szego}, p. 59]
\begin{equation}\label{symmetry1}
P_k^{(\alpha,\beta)}(-x)=(-1)^k P_k^{(\beta,\alpha)}(x).
\end{equation}
Let $\mathbf D^{(\alpha,\beta)}=[d_{i,j}]_{i,j\geq 0}$ and $\mathbf D^{(\beta,\alpha)}=[\widetilde d_{i,j}]_{i,j\geq 0}$
be the matrix representation of $-\mathcal L_c^{(\alpha,\beta)}$ and $-\mathcal L_c^{(\beta,\alpha)}$ with respect to the 
basis of Jacobi polynomials $\wJ_k$ and $\widetilde P_k^{(\beta,\alpha)},$ respectively. Then from \eqref{recursion1} and \eqref{eigensystem2}, one gets 

\begin{equation}\label{coefficientsD}
\widetilde d_{i,i+k}= (-1)^ d_{i,i+k},\quad  -2\leq k\leq 2,\quad i\geq 0\mbox{ and } \widetilde d_{i,j}=0,\quad \mbox{ if \ \ } |i-j|\geq 3.
\end{equation}
Let $\mathbf B_n = [\beta_k^n, k \geq 0]^T$ and $\widetilde{\mathbf B}_n = [\widetilde \beta_k^n, k\geq 0]^T =[(-1)^k \beta_k^n, k\geq 0]^T. $
Since $\mathbf D^{(\alpha,\beta)} \mathbf B_n = \chi_n(c) \mathbf B_n,$ then by using \eqref{coefficientsD}, it is easy to see that 
$\mathbf D^{(\beta,\alpha)} \widetilde{\mathbf B}_n = \chi_n(c) \mathbf B_n.$ This means that $\psi_{n,c}^{(\alpha,\beta)}$ and 
$\psi_{n,c}^{(\beta,\alpha)}$ are associated to the same eigenvalue $\chi_n(c).$ Moreover, the series expansion of $\psi_{n,c}^{(\beta,\alpha)}$
in the basis $ \{\widetilde P_k^{(\beta,\alpha)},\,\, k\geq 0\}$ is obtained from the series expansion of $\psi_{n,c}^{(\alpha,\beta)},$ as follows
\begin{equation}
\label{series2}
\psi_{n,c}^{(\beta,\alpha)}(x)= c_n \sum_{k=0}^{\infty} (-1)^k \beta_k^n \widetilde P_k^{(\beta,\alpha)}(x),
\end{equation}
for some constant $c_n.$ 
By combining \eqref{symmetry1} and the previous equality, one concludes that $\psi_{n,c}^{(\alpha,\beta)}(-x)= c_n \psi_{n,c}^{(\beta,\alpha)}(x).$
Moreover, since $\|\psi_{n,c}^{(\alpha,\beta)}\|_{L^2(I,\omega_{\alpha,\beta})}=\|\psi_{n,c}^{(\beta,\alpha)}\|_{L^2(I,\omega_{\beta,\alpha})}=1$
and since $\psi_{n,c}^{(\alpha,\alpha)}$ has the same parity as $n,$ see \cite{Wang2},  then $c_n=(-1)^n.$ This concludes the proof
of \eqref{symmetry}.$\qquad \Box$

Also, we should mention that the $(n+1)-$th eigenvalue $\chi_n(c)$ satisfies the following classical inequalities,
\begin{equation}
\label{boundschi}
n (n+\alpha+\beta+1) \leq \chi_n(c) \leq n (n+\alpha+\beta+1) +c^2,\quad \forall n\geq 0.
\end{equation}

To get the previous upper  bound, we  consider
the following Sturm-Liouville form of $-\mathcal L_c^{(\alpha,\beta)},$
\begin{equation}\label{Sturm_form}
- \mathcal L_c^{(\alpha,\beta)} (f)(x)= - \frac{d}{dx}\left[ \omega_{\alpha,\beta}(x) (1-x^2) f'(x)\right] +c^2 x^2 \omega_{\alpha,\beta}(x)= 
- \mathcal L_0^{(\alpha,\beta)} (f)(x) +c^2 x^2 \omega_{\alpha,\beta}(x).
\end{equation}
Then, from the well known Poincar\'e Min-Max characterization of the eigenvalue of a  self-adjoint  operator, 
applied to the operator $-\mathcal L_c^{(\alpha,\beta)}, $ one gets 
\begin{eqnarray*}
\label{bound1}
\chi_n(c)&=& \min_{\mbox{ dim } H =n} \,\,  \max_{u\in H, \|u\|=1} \int_{-1}^1 \left(-\mathcal L_0^{(\alpha,\beta)}(u)(x) u(x) +c^2 x^2 (u(x))^2 \right)
\omega_{\alpha,\beta}(x)\, dx \nonumber\\
&\leq & \min_{ \mbox{ dim } H =n} \,\,  \max_{u\in H, \|u\|=1} \int_{-1}^1 -\mathcal L_0^{(\alpha,\beta)}(u)(x) u(x) \omega_{\alpha,\beta}(x)\, dx +c^2 \| u\|^2\\
&\leq & \chi_n(0)+ c^2= n(n+\alpha+\beta+1)+c^2.
\end{eqnarray*}
Next, to get a lower bound, it suffices to see that the self-adjoint operator $-\mathcal L_c^{(\alpha,\beta)}-(-\mathcal L_0^{(\alpha,\beta)})= c^2 x^2$ is a positive operator, which implies that $\chi_n(c)\geq \chi_n(0).$

\subsection{Local estimates of the eigenfunctions of $\mathcal L_c^{(\alpha,\beta)}.$}

In this paragraph, we give various explicit and local estimates of the $\psi_{n,c}^{(\alpha,\beta)}.$ These estimates will be needed to prove some of the results of sections 3 and 4 of this work. 
We should mention that in the literature, only few references have studied the  problem of the 
explicit estimates of the classical PSWFs and their eigenvalues $\chi_n(c),$ 
 see \cite{Bonami-Karoui1, Bonami-Karoui2, Osipov}. 
The following proposition provides us with explicit local bounds
of $\psi_{n,c}^{(\alpha,\beta)},$ generalizing a similar  result given 
in  \cite{Bonami-Karoui1} for the special case  $\alpha=\beta=0.$

\begin{proposition}
For real numbers $c>0,\, \alpha, \beta > -1,$ with  $\alpha+\beta+1 \geq 0.$  Let  $n\in \mathbb N$  be such $q= c^2/\chi_n(c) <1.$ Then we have
\begin{equation}\label{localest1}
\sup_{t\in [0,1]} (1-t^2) \omega_{\alpha,\beta}(t) \left(|\psi_{n,c}^{(\alpha,\beta)}(t)|^2+ \frac{1-t^2}{(1-q t^2)\chi_n(c)}|(\psi_{n,c}^{(\alpha,\beta)})'(t)|^2\right) \leq 2(1+\max(\alpha,\beta)).
\end{equation}
Moreover, if $\alpha=\beta,$ then we have 
\begin{equation}\label{localest2}
\sup_{t\in [-1,1]} (1-t^2) \omega_{\alpha}(t) \left(|\psi_{n,c}^{(\alpha)}(t)|^2+ \frac{1-t^2}{(1-q t^2)\chi_n(c)}|(\psi_{n,c}^{(\alpha)})'(t)|^2\right) \leq 1+\alpha.
\end{equation}
\end{proposition}

\noindent
{\bf Proof:} The proof uses a classical technique for the local estimates of the eigenfunctions
of a Sturm-Liouville operator. In our case, we first note that by using property $(P_2)$ of proposition \ref{properties1}, it suffices to consider the case  $\alpha \geq \beta,$
since the case $\beta  \geq \alpha,$ follows from the equality \eqref{symmetry}.
Next,  consider the auxiliary function, defined on $[0,1]$ by 
$$ Z_n(t)= \left(\psi_{n,c}^{(\alpha,\beta)}(t)\right)^2+ \frac{1-t^2}{\chi_n(c) (1-q t^2)}\left((\psi_{n,c}^{(\alpha,\beta)})'\right)^2(t).$$
Since $ \psi_{n,c}^{(\alpha,\beta)}$ is the eigenfunction of the operator $-\mathcal L^{(\alpha,\beta)}_c$ associated with the eigenvalue
$\chi_n(c),$ then  straightforward computations give us
\begin{equation}
\label{deriv1}
 Z_n'(t)= \frac{2 ((\psi_{n,c}^{(\alpha,\beta)})')^2(t)}{\chi_n(c) (1-q t^2)}\left( (\alpha-\beta)+(1+\alpha+\beta) t + qt
\frac{1-t^2}{1-q t^2}\right).
\end{equation}
Since $0\leq q\leq 1,$ $\alpha-\beta\geq 0$ and $\alpha+\beta+1\geq 0,$ then it is easy to see that 
$$Z_n'(t)\geq 0,\quad \forall\,\,  t\in [0,1].$$
Next, we consider a second auxiliary function, given by 
$$ K_n(t)= (1-t^2) \omega_{\alpha,\beta}(t) Z_n(t),\quad t\in [0,1].$$ Then, by using \eqref{deriv1}, one can easily check that there exists 
a positive valued function $A(\cdot)$ on $[-1,1]$ with 
\begin{eqnarray*}
K_n'(t)&=& \omega_{\alpha,\beta}(t)\left( -2 t +(\beta-\alpha)-(\alpha+\beta)t\right) Z_n(t)+ (1-t^2)\omega_{\alpha,\beta}(t) Z'_n(t)\\
&= & \left( -2 t +(\beta-\alpha)-(\alpha+\beta)t\right)\omega_{\alpha,\beta}(t) \left(\psi_{n,c}^{(\alpha,\beta)}\right)^2(t)+ A(t)\left(\psi_{n,c}^{(\alpha,\beta)})'\right)^2(t)\\
&\geq &  \left((\beta-\alpha)-(2+\alpha+\beta)t\right)\omega_{\alpha,\beta}(t) \left(\psi_{n,c}^{(\alpha,\beta)}\right)^2(t).
\end{eqnarray*}
Finally, since $K_n(1)=0$ and since ${\displaystyle \int_{-1}^1 \left(\psi_{n,c}^{(\alpha,\beta)}\right)^2(t) \, \omega_{\alpha,\beta}(t)\, dt =1,}$ then 
by using the last inequality, one gets 
$$K_n(t)-K_n(1)=K_n(t) \leq \max_{t\in [0,1]} \left((\alpha-\beta)+(2+\alpha+\beta)t\right)
\int_{0}^1 \left(\psi_{n,c}^{(\alpha,\beta)}\right)^2(t) \, \omega_{\alpha,\beta}(t)\, dt \leq 2(1+\alpha).$$
Finally, if $\beta=\alpha,$ then from the parity of $\psi_{n,c}^{(\alpha)}(t),$ we have 
${\displaystyle \int_{0}^1 \left(\psi_{n,c}^{(\alpha)}\right)^2(t) \, \omega_{\alpha}(t)\, dt=1/2,}$ which means that the previous
upper bound is replaced by $1+\alpha.$ This concludes the proof of the proposition.

The following proposition provides us with an estimate of the maximum of the $\psi_{n,c}^{(\alpha,\beta)}$ inside the interval $I.$

\begin{proposition} Let $c>0,$ and $\alpha\geq  \beta$ with $\alpha+\beta \geq -1,$ then for any positive integer $n$ with $q=c^2/\chi_n(c)\leq 1,$
we have 
\begin{equation}\label{maxpsi1}
\sup_{x\in [0,1]} |\psi_{n,c}^{(\alpha,\beta)}(x)|=|\psi_{n,c}^{(\alpha,\beta)}(1)|\leq  C_{\alpha} (\chi_n(c))^{\frac{1+\alpha}{2}},\quad C_{\alpha}=
\frac{2^{1+\max(\alpha,0)}}{\sqrt{3+\alpha}}\Big(\frac{1+\alpha}{3+\alpha} \Big)^{1+\frac{\alpha}{2}}
\end{equation}
Moreover, if $\alpha =\beta,$ then we have 
\begin{equation}\label{maxpsi2}
\sup_{x\in [-1,1]} |\psi_{n,c}^{(\alpha)}(x)|=|\psi_{n,c}^{(\alpha)}(1)|\leq  \frac{C_{\alpha}}{\sqrt{2}} (\chi_n(c))^{\frac{1+\alpha}{2}}.
\end{equation}
\end{proposition}

\noindent
{\bf Proof:}  
We first recall  that the auxiliary function $Z_n$ given by :
$$ Z_n(x)= \left(\psi_{n,c}^{(\alpha,\beta)}(x)\right)^2+ \frac{1-x^2}{\chi_n(c) (1-q x^2)}\left((\psi_{n,c}^{(\alpha,\beta)})'\right)^2(x)$$ 
is increasing over [0,1] whenever $q=c^2/\chi_n(c)\leq 1,$ $\alpha \geq \beta$ and $\alpha+\beta+1 \geq 0.$
Hence, we  have $${\displaystyle \sup_{x\in [0,1]} Z_n(x)= Z_n(1)=  \left(\psi_{n,c}^{(\alpha,\beta)}(x)\right)^2}$$ which implies that 
\begin{equation}\label{eq2}
\sup_{x\in [0,1]}|\psi_{n,c}^{(\alpha,\beta)}(x)|=|\psi_{n,c}^{(\alpha,\beta)}(1)|,\quad \forall n\in \mathbb N, \mbox{ with } q\leq 1.
\end{equation}
Moreover, if $\alpha=\beta$  then from the parity of $\psi_{n,c}^{(\alpha)},$ one gets
\begin{equation}\label{eq3}
\sup_{x\in [-1,1]}|\psi_{n,c}^{(\alpha)}(x)|=|\psi_{n,c}^{(\alpha)}(1)|,\quad \forall n\in \mathbb N, \mbox{ with } q\leq 1.
\end{equation}
Next, we show how to get the upper bounds of $ |\psi_{n,c}^{(\alpha,\beta)} (1)| $ and $ |\psi_{n,c}^{(\alpha)} (1)|. $
To alleviate notation, we simply denote $\psi_{n,c}^{(\alpha,\beta)}$ by $\psi_{n,c}$ and $\chi_n(c)$ by $\chi_n.$ Also, without loss of generality, we may assume that $\psi_{n,c}(1)>0.$
Since
\begin{equation}
(\psi_{n,c}'(x) (1-x^2) \omega_{\alpha,\beta}(x))'=-\chi_n(c) \omega_{\alpha,\beta}(x)(1-qx^2) \psi_{n,c} (x),
\end{equation}
then 
\begin{eqnarray*}
\psi_{n,c}'(x) &=& \frac{\chi_n}{(1-x^2) \omega_{\alpha,\beta}(x)} \int_x^1 \omega_{\alpha,\beta}(t)(1-qt^2) \psi_{n,c} (t) dt\\
& \leq&  \frac{\chi_n}{(1-x^2)}(1-qx^2)(1-x) \psi_{n,c}(1) =  \chi_n(1-qx^2) \psi_{n,c}(1).
\end{eqnarray*}
Hence,  
\begin{equation}\label{Equat1}
\psi_{n,c}(1)-\psi_{n,c}(x) \leq \chi_n Q_q(x) \psi_{n,c}(1),\qquad Q_q(x)=(1-qx^2)(1-x^2).
\end{equation}
Next, let $x_n\in [0,1]$ be such that $ Q_q(x_n)= \frac{a}{\chi_n},$ where the constant $a$ to be fixed later on.
By substituting $x$ with $x_n$ in \eqref{Equat1} and by using \eqref{localest1}, one gets 
$$\psi_{n,c}(1) \leq  \frac{1}{1-a} \psi_{n,c}(x_n) \leq \frac{1}{(1-a)}\sqrt{2(1+\alpha)}\left(\frac{1}{Q_q(x_n)\omega_{\alpha,\beta}(x_n)}\right)^{1/2}.$$
That is 
\begin{equation}\label{Equat2}
\psi_{n,c}(1) \leq \frac{\sqrt{2(1+\alpha)}}{a^{1/2}(1-a)}\frac{\chi_n^{1/2}}{\sqrt{\omega_{\alpha,\beta}(x_n)}} 
\end{equation}
Note that  the admissible solution of  $ Q_q(x_n)=\frac{a}{\chi_n}$  is given by 
$ x_n=\left( \frac{(q+1)-\sqrt{(q-1)^2+\frac{4aq}{\chi_n}}}{2q} \right)^{1/2}.$
Consequently,
$$ \frac{1-\frac{a}{\chi_n}}{1+\sqrt{\frac{aq}{\chi_n}}}\leq x_n = \left( \frac{2(1-\frac{a}{\chi_n})}{q+1+(1-q)\sqrt{1+\frac{4aq}{\chi_n(q-1)^2}}} \right)^{\frac{1}{2}}\leq \Big(1-\frac{a}{\chi_n}\Big)^{\frac{1}{2}}\leq 1-\frac{a}{2\chi_n}.$$
It is easy to see that in this case, we have 
$$\frac{a}{2\chi_n} \leq 1-x_n \leq \frac{a}{\chi_n}\left(1+\sqrt{\frac{\chi_n}{a}}\right).$$
Consequently, by using  the first inequality when $\alpha\geq 0$ and the second inequality when 
$-1/2\leq \alpha<0,$ one gets 
\begin{equation}\label{Equat3}
\frac{1}{\sqrt{\omega_{\alpha,\beta}(x_n)}} \leq \left\{\begin{array}{ll} \left(\frac{\chi_n}{2 a}\right)^{\alpha/2} 2^{-\beta/2}
\leq \sqrt{2} \left(\frac{\chi_n}{a}\right)^{\alpha/2} &\mbox{ if } -1/2\leq \alpha<0\\
& \\
\left(\frac{2\chi_n}{a}\right)^{\alpha/2} \frac{1}{(1+x_n)^{\beta/2}} \leq 2^{\alpha+1/2}\left(\frac{\chi_n}{a}\right)^{\alpha/2}&\mbox{ if }
\alpha \geq 0.\end{array}\right.
\end{equation}
Hence, by combining \eqref{Equat2} and \eqref{Equat3}, one gets
\begin{equation}
\psi_{n,c}^{(\alpha,\beta)}(1) \leq \frac{2^{1+\max(\alpha,0)}\sqrt{(1+\alpha)}}{a^{(1+\alpha)/2}(1-a)}\chi_n^{\frac{1+\alpha}{2}}.
\end{equation}
Since the  maximum of $ a^{\gamma}(1-a) $ is attained at $ a= \frac{\gamma}{1+\gamma} $
then for $\gamma=\frac{1+\alpha}{2},$ one gets \eqref{maxpsi1}. Finally, \eqref{maxpsi2} follows from the parity of $\psi_{n,c}^{(\alpha)}$
and \eqref{localest2}.$\qquad \Box$\\

\noindent
\begin{remark}
The techniques used for the proof of  inequality \eqref{Equat2} are similar to those used in \cite{Bonami-Karoui1}
to prove a similar inequality, restricted to the special case $\alpha=\beta=0.$
Nonetheless, the general setting of the previous proposition requires handling  the new quantity 
$\sqrt{\omega_{\alpha,\beta}(x_n)}$ that generates extra difficulties to obtain the local estimates \eqref{maxpsi1}
and \eqref{maxpsi2}.
\end{remark}

\section{Generalized prolate spheroidal wave functions: Computations and analytic extension.}

In the sequel, we restrict ourselves to the case $\alpha=\beta> -1.$ In the first part of this section, we further 
improve  the  super-exponential decay rate  
of  the GPSWFs expansion coefficients $(\beta_k^n)_k,$ that  has been given in given in \cite{Wang2}.
Then, we show that for sufficiently large values of $n$ and
 up a certain order $K_n,$ all the coefficients $\beta_k^n,\, 0\leq k\leq K_n$ are positive. As a consequence of this positivity result
 and the previous fast decay of the $\beta_k^n,$ we show that in the case where $\alpha=\beta,$ the expansion coefficients $(\beta_k^n)_k$
 are  essentially concentrated around $k=n.$ In the second part of this section, we give the analytic extension of the GPSWFs, together with an
 explicit expression for the eigenvalues $\mu_n^{(\alpha)}(c)$ as a ratio of two fast  convergent series.

\subsection{Computation and analytic extension of the GPSWFs}

 We first  note that in the interesting special case where $\alpha=\beta,$ formula \eqref{fourier_Jacobi} is simplified in a significant manner. This is given 
 by the following lemma.
 \begin{proposition} Let $\alpha > -1,$ then we have 
 \begin{equation}\label{fourier_Jacobi2}
 \int_{-1}^{1} e^{ixy} P_k^{(\alpha , \alpha)}(y) \omega_{\alpha}(y)\,  dy= i^k \sqrt{\pi} \left(\frac{2}{x}\right)^{\alpha+1/2}
 \frac{\Gamma(k+\alpha+1)}{\Gamma(k+1)}J_{k+\alpha+1/2}(x),\quad x\in \mathbb R.
 \end{equation}
 Here, $J_{a}$ denotes the Bessel function of the first kind and order $\alpha.$
 \end{proposition}
 
 \noindent
 {\bf Proof:} It is well known, see for example [\cite{Andrews}, p. 200], that if 
  $a > -1$ and $z= - 2 i  x,\,\, x\in \mathbb R,$ then we have,  
  \begin{equation}\label{Kummer2}
 _1F_1(a +1/2,2 a +1; -2ix)= \Gamma(a+1)\left(\frac{2}{x}\right)^{a}  e^{-ix} J_{a}(x).
 \end{equation}
 By combining \eqref{fourier_Jacobi} and the previous equality with $a=k+\alpha+1/2,$ one gets
 \begin{eqnarray}\label{fourier_Jacobi3}
 \mathcal F_1^{\alpha}(P_k^{(\alpha , \alpha)})(x)&=&\int_{-1}^{1} e^{ixy} P_k^{(\alpha , \alpha)}(y) \omega_{\alpha}(y)\,  dy\nonumber\\
 &=& \frac{(i x)^k2^{2k+3\alpha+3/2} }{x^{k+\alpha+1/2}}
  \frac{\Gamma(k+\alpha+3/2)}{\Gamma(k+1)} B(k+\alpha+1,k+\alpha+1) J_{k+\alpha+1/2}(x).
 \end{eqnarray}
 Moreover, by using the following identities of Beta and Gamma functions,
 $$ B(x,y)=\frac{\Gamma(x) \Gamma(y)}{\Gamma(x+y)},\quad \Gamma(x+1)=x\Gamma(x),\quad \Gamma(b)\Gamma(b+1/2)=\sqrt{\pi} \frac{\Gamma(2b)}{2^{2b-1}},$$
 one gets
 \begin{eqnarray*}
 \mathcal F_1^{\alpha} P_k^{(\alpha , \alpha)}(x)&=& i^k \frac{2^{2k+3\alpha+3/2}}{(x)^{\alpha+1/2}}  (k+\alpha+1/2) \frac{\Gamma(k+\alpha+1)}{\Gamma(k+1)}\frac{\Gamma(k+\alpha+1)\Gamma(k+\alpha+1/2)}{\Gamma(2k+2\alpha+2)}J_{k+\alpha+1/2}(x)\\
 &=& i^k\frac{\Gamma(k+\alpha+1)}{\Gamma(k+1)}\frac{2^{2k+3\alpha+3/2}}{(cx)^{\alpha+1/2}}\frac{k+\alpha+1/2}{2k+2\alpha+1}
 \frac{\Gamma(k+\alpha+1)\Gamma(k+\alpha+1/2)}{\Gamma(2k+2\alpha+1)}J_{k+\alpha+1/2}(x)\\
 &=&i^k \sqrt{\pi} \left(\frac{2}{x}\right)^{\alpha+1/2}
 \frac{\Gamma(k+\alpha+1)}{\Gamma(k+1)}J_{k+\alpha+1/2}( x),\quad x\in \mathbb R.\qquad\qquad \Box
 \end{eqnarray*}
 
 \begin{remark}
 In the special case $\alpha=0,$ the equality \eqref{fourier_Jacobi2} is reduced to the well known classical finite Fourier transform
 of Legendre function, see for example [\cite{Andrews}, p. 343].
 \end{remark}
 
 As a first consequence of the previous proposition, one gets a simple and straightforward proof of the following 
 result that has been already  given by  Lemma 3.4 in \cite{Wang2} and with different kind of proof.

 \begin{corollary} Under the above notations, for any real numbers $c>0$ and $\alpha>-1,$ we have
 \begin{equation} \label{beta2}
 \beta_k^n= \int_{-1}^1 \widetilde P_k^{(\alpha,\alpha)}(y) \psi_{n,c}^{(\alpha)}(y) \omega_{\alpha}(y)\, dy=
 \frac{2\sqrt{\pi} i^k}{\mu_n^{(\alpha)}(c) \sqrt{h_k}}\left(1+(-1)^{n+k}\right) \int_0^1 \frac{J_{k+\alpha+\frac{1}{2}}(ct)}{(ct)^{\alpha+\frac{1}{2}}} \psi_{n,c}^{(\alpha)}(t) \omega_{\alpha}(t)\, dt.
 \end{equation}
 \end{corollary}
 
 \noindent
 {\bf Proof:} Just write
 \begin{eqnarray*}
 \beta_k^n &=& \frac{1}{\mu_n^{(\alpha)}(c)}\int_{-1}^1 \widetilde P_k^{(\alpha,\alpha)}(y) \left(\int_{-1}^1 e^{icty} \psi_{n,c}^{(\alpha)}(t) \omega_{\alpha}(t)\, dt \right) \omega_{\alpha}(y)\, dy\\
 &=& \frac{\sqrt{\pi} i^k}{\mu_n^{(\alpha)}(c) \sqrt{h_k}} \int_{-1}^1  \frac{J_{k+\alpha+\frac{1}{2}}(ct)}{(ct)^{\alpha+\frac{1}{2}}} \psi_{n,c}^{(\alpha)}(t) \omega_{\alpha}(t)\, dt.
 \end{eqnarray*}
 To conclude, it suffices to write the previous  integral as $\int_{-1}^1 = \int_{0}^1 + \int_{-1}^0$ and use  the facts that the function $\psi_{n,c}^{(\alpha)}$ and ${\displaystyle t\mapsto \frac{J_{k+\alpha+\frac{1}{2}}(ct)}{(ct)^{\alpha+\frac{1}{2}}}}$ has the same parity as $n$ and $k,$ respectively.$\qquad \Box$

  A decay rate of the expansion coefficients is given by the following proposition that improves the result given by 
  Theorem 3.4 in \cite{Wang2}.
 
 \begin{proposition}\label{prop_decay1}
 For given  real numbers $c>0,\, \,  \alpha >- 1$ and integers  $k,n \in \mathbb N,$ let  $${\displaystyle \beta_{k}^{n}=\int_{-1}^{1} \widetilde P_k^{(\alpha,\alpha)}(x) \; \psi^{(\alpha)}_{n,c}(x)\, \omega_{\alpha}(x)  dx.}$$ Then, we have 
 \begin{equation}\label{decay1_coeff} 
 |\beta_k^n|\leq\frac{C_{\alpha}}{|\mu_n^{(\alpha)}(c)|}\frac{1}{k^{1+\alpha/2} }\frac{1}{2^k} \left(\frac{ec}{2k+1}\right)^k,
  \end{equation}
 where ${\displaystyle  C_{\alpha}=\frac{\pi^{7/4} \sqrt{\Gamma(1+\alpha)}(3/2)^{3/4}(3/2+2\alpha)^{3/4+\alpha}}{2^{\alpha+1}e^{\alpha+5/4}}.}$
 \end{proposition}

\noindent
{\bf Proof:}  By using the expression of $\beta_k^n$ and by combining 
(\ref{eq1}) and (\ref{fourier_Jacobi}), one gets   
\begin{eqnarray*}
\beta_k^n &=&\frac{1}{\mu_n^{(\alpha)}(c)}\int_{-1}^{1}\int_{-1}^{1}\wJ_{k}(x) \left(\int_{-1}^1 e^{icxy} \psi_{n,c}(y) \omega_{\alpha}(y) \, dy\right)  \omega_{\alpha}(x)\, dx\\
&=& \frac{1}{\mu_n(c)}\int_{-1}^{1}\left(\int_{-1}^{1} e^{icxy}\wJ_{k}(x) \omega_{\alpha}(x)\, dx\right)
\psi_{n,c}^{(\alpha)}(y) \omega_{\alpha}(y)\, dy\\
&=& \frac{1}{\mu_n^{(\alpha)}(c)}\int_{-1}^{1}
\frac{(icy)^ke^{icy}}{k! \sqrt{h_k}}2^{k+2\alpha+1}B(k+\alpha+1,k+\alpha+1)\cdot\\
&&\qquad\qquad\qquad _1F_1(k+\alpha+1,2k+2\alpha+2;-2icy)\psi_n(y)\, \omega_{\alpha}(y)\, dy.
\end{eqnarray*}
On the other hand, from the integral representation of Kummer's function given by \eqref{Kummer1}, one gets 
\begin{eqnarray*}
\lefteqn{|_1F_1(k+\alpha+1,2k+2\alpha+2;-2icy)|=\frac{\Gamma(2k+2\alpha+2)}{\Gamma(k+\alpha+1)^2}
\left|\int_0^1 e^{-2icyt}t^{k+\alpha}(1-t)^{k+\alpha}dt\right|}\\
&& \qquad\qquad \leq  \frac{1}{B(k+\alpha+1,k+\alpha+1)}\left|\int_0^1 e^{-2icyt}t^{k+\alpha}(1-t)^{k+\alpha}dt\right|\leq \frac{B(k+\alpha+1,k+\alpha+1)}{B(k+\alpha+1,k+\alpha+1)}=1.
\end{eqnarray*}
Consequently, we have 
\begin{eqnarray}\label{Ineq2.1}
|\beta_k^n| &\leq& \frac{B(k+\alpha+1,k+\alpha+1)}{|\mu_n^{(\alpha)}(c)|}\frac{c^k 2^{k+2\alpha+1}}{k! \sqrt{h_k}}
\left(\int_{-1}^1  (\psi_{n,c}^{(\alpha)})^2(y) \omega_{\alpha}(y)\, dy \right)^{1/2} \left(\int_{-1}^1 
y^{2k}\omega_{\alpha}(y)\, dy\right)^{1/2} \nonumber\\
&\leq & \frac{B(k+\alpha+1,k+\alpha+1)}{|\mu_n^{(\alpha)}(c)|}\frac{c^k 2^{k+3/2(\alpha+\alpha+1)}}{k! \sqrt{h_k}} \sqrt{I_{2k}^{\alpha,\alpha}}.
\end{eqnarray}
Note that from \eqref{Ineq2}, we have 
\begin{equation}\label{Ineq3}
\frac{c^k}{k!}=\frac{c^k}{\Gamma(k+1)}\leq \frac{1}{\sqrt{2}}\frac{1}{\sqrt{2k+1}} \left(\frac{ec}{2k+1}\right)^k. 
\end{equation}
In a similar manner, we get the following upper bound and lower bound of the quantity $B(k+\alpha+1,k+\alpha+1)$ and the normalization
constant $h_k,$ given as follows.
\begin{eqnarray}\label{Ineq6}
B(k+\alpha+1,k+\alpha+1)&=&\frac{\Gamma(k+\alpha+1) ^2}{\Gamma(2k+2\alpha+2)} 
\leq  \frac{\sqrt{2}\pi}{2^{2k+2\alpha+1}}\frac{(2k+2\alpha+1)^{2k+2\alpha+1}}{(2k+2\alpha+3/2)^{2k+2\alpha+3/2}} \nonumber\\
&\leq & \frac{\sqrt{2}\pi}{2^{2k+2\alpha+1}} \frac{1}{\sqrt{2k+2\alpha+3/2}}
\end{eqnarray}
Also, by using \eqref{Ineq2}, the decay of the function $\varphi,$ given by \eqref{function} as well some  straightforward computations, one gets 
\begin{equation}\label{Ineq7}
\frac{e}{\pi} \frac{(2e)^{2\alpha+1}}{(3/2)^{3/2}(3/2+2\alpha)^{3/2+2\alpha}}\frac{1}{2k+2\alpha+1}\leq h_k \leq
\frac{\pi}{e}\left(\frac{2}{e}\right)^{2\alpha+1}\frac{(3/2)^{3/2}(3/2+2\alpha)^{3/2+2\alpha}}{2k+2\alpha+1}.
\end{equation}
Finally, by combining \eqref{Ineqq5}, \eqref{Ineq2.1}, \eqref{Ineq3}--\eqref{Ineq7}, one gets the desired result \eqref{decay1_coeff}.$\qquad \Box$

\begin{remark} By using our notation, the decay rate of the $(\beta_k^n)_k,$ given by Theorem 3.4 of \cite{Wang2} can be written as 
${\displaystyle \frac{C^{"}_{\alpha}}{|\mu_n^{(\alpha)}(c)|}\frac{1}{k^{1+\alpha/2} } \left(\frac{ec}{2k+1}\right)^k,}$ for some constant $C^{"}_{\alpha}.$
The previous proposition ensures that this decay is further improved by a factor of $1/2^k.$
\end{remark}

The following theorem provides us with a second decay rate of the $(\beta_k^n)_{k\geq 0},$ valid for sufficiently large values of $n$ and 
 the values of $0\leq k < n$  not too close to $n.$ We should mention that the  techniques of  the proof of this theorem, given in the the Appendix, are inspired from those developed for the special case $\alpha=0$ and given in a joint work of one of us \cite{Bonami-Karoui4}.

\begin{theorem} Let $c>0,$ be a fixed positive real number. Then, for all positive integers  $n,k$ such that 
$q=c^2/\chi_n \leq 1$ and $k(k+2\alpha+1)+C_{\alpha}\,  c^2\leq \chi_n(c)$, we have
\begin{equation}
\label{Decay2beta}
|\beta_0^n|\leq \sqrt{\frac{\Gamma(\alpha+3/2)}{\sqrt{\pi}\Gamma(\alpha+1)}} \sqrt{1+\alpha} |\mu_n^{(\alpha)}(c)|\quad \mbox{ and }\quad |\beta_k^n|\leq C'_{\alpha} \Big( \frac{2}{q}\Big)^k|\mu_n^{(\alpha)}(c)|.
\end{equation}
Here,  ${\displaystyle C'_{\alpha} =\frac{2^{\alpha}(3/2)^{3/4}(3/2+2\alpha)^{3/4+\alpha}}{e^{2\alpha+3/2}}\sqrt{1+\alpha}}$  and $C_{\alpha}= 2 M_{\alpha}+N_{\alpha}$ with 
\begin{equation}\label{C_alpha}
 M_{\alpha}=\max\Big({1/4,\sqrt{\frac{2(2\alpha+2)}{(2\alpha+5)(2\alpha+3)^2}}}\Big)
 ,\,\, N_{\alpha}=\max \Big( \frac{3}{2\alpha+5}, \frac{1}{2}+\frac{|4\alpha^2-1|}{(2\alpha+3)(2\alpha+7)}\Big).
 \end{equation}
 \end{theorem}

\subsection{Analytic extension of the GPSWFs}

In this paragraph, we give explicit formulae for the analytic extension of the GPSWFs to the whole real line, as well for computing
the eigenvalues $\mu_n^{(\alpha)}(c)$ associated with the weighted finite Fourier transform $\mathcal F_c^{\alpha}.$ We first note that due to 
equation \eqref{eq1}, the GPSWFs have analytic extension to $\mathbb R.$ In fact, it is well known, see 
[\cite{Szego}, p.168], that 
if $\alpha > -1,$ then 
$${\displaystyle \sup_{y\in [-1,1]} |\wJ_k(y)|= \leq 
\frac{1}{\sqrt{h_k} B(k+1,\alpha+1)}\leq M_{\alpha} k^{\alpha+1/2},}$$ 
for some constant $M_{\alpha}.$ Moreover, by using   the super-exponential decay rate of the expansion coefficients $(\beta_k^n)_k$
combined with \eqref{eq1}, \eqref{expansion1} and \eqref{fourier_Jacobi2}, one gets 
\begin{eqnarray*}
\psi_{n,c}^{(\alpha)}(x)&=&\frac{\sqrt{\pi} 2^{\alpha+1/2}}{\mu_n^{(\alpha)}(c)} \sum_{k\geq 0} i^k \beta_k^n\, \frac{\Gamma(k+\alpha+1)}{\sqrt{h_k} k!}
\frac{J_{k+\alpha+1/2}(c x)}{(cx)^{\alpha+1/2}},\quad \forall\, x\neq 0.
\end{eqnarray*}
Moreover, from the parity of the $\psi_{n,c}^{(\alpha)},$   it is easy to see that $\mu_n^{(\alpha)}(c)= i^n |\mu_n^{(\alpha)}(c)|,$ $n\geq 0.$
Hence, by using the fact that  the previous expansion coincides with the expansion \eqref{expansion1} at $x=1,$ one obtains the following analytic extension of the GPSWFs as well as an explicit formula for  their associated eigenvalues $\mu_n^{(\alpha)}(c),$
\begin{equation}
\label{expansion2}
\psi_{n,c}^{(\alpha)}(x)=\frac{\sqrt{\pi} 2^{\alpha+1/2}}{|\mu_n^{(\alpha)}(c)|} \sum_{k\geq 0} i^{k-n} \beta_k^n\, \frac{\Gamma(k+\alpha+1)}{\sqrt{h_k} k!}
\frac{J_{k+\alpha+1/2}(c x)}{(cx)^{\alpha+1/2}},\quad \forall\, x\neq 0
\end{equation}
with
\begin{equation}
\label{mu_n}
\mu_n^{(\alpha)}(c)= i^n \sqrt{\pi}\left(\frac{2}{c}\right)^{\alpha+1/2}\frac{\sum_{k\geq 0} i^{k-n} \beta_k^n\, \frac{\Gamma(k+\alpha+1)}{\sqrt{h_k} k!}
{J_{k+\alpha+1/2}(c)}}{\sum_{k\geq 0}  \frac{\beta_k^n}{\sqrt{h_k}\, B(k+\alpha,k) }},\quad n\geq 0.
\end{equation}
We should mention that due to the facts that the coefficients $(\beta_k^n)_k$ are concentrated around $k=n$ and decay super-exponentially,
 the previous formula is  accurate  and  practical  for computing the $\mu_n^{(\alpha)}(c).$  Also,
note that in \cite{Wang2}, the authors have given some properties of the eigenvalues $\mu_n^{(\alpha)}(c)$ (denoted by  $\lambda_n^{(\alpha)}(c)$ in \cite{Wang2}).
In particular, by considering the operator $\mathcal F_c^{\alpha}$ as a Hilbert-Schmidt operator acting on $L^2(I, \omega_{\alpha}),$
it has been shown that 
$$\sum_{n\geq 0} |\mu_n^{(\alpha)}(c)|^2=\| \mathcal F_c^{\alpha}\|_{HS}^2 = \left(\int_{-1}^1 \omega_{\alpha}(x)\, dx \right)^2= \frac{\pi \Gamma^2(1+\alpha)}{\Gamma^2(\alpha+3/2)}.$$
Here, $\|\cdot\|_{HS}$ denotes the Hilbert-Schmidt norm. More importantly, in \cite{Wang2}, 
the authors have noted   that the $\mu_n^{(\alpha)}(c)$ has an asymptotic super-exponential decay rate given by 
\begin{equation}
\label{decaymu_n}
|\mu_n^{(\alpha)}(c)| \approx \frac{e^{\alpha}}{4^{\alpha}}\sqrt{\frac{\pi e}{2n+2\alpha+3}}\left(\frac{ec}{4n+4\alpha+2}\right)^n,\quad n>>1.
\end{equation}

\section{GPSWFs as solutions of an energy maximization problem and  quality of approximation.}

In the first part of this section, we show that in the case where $\alpha \geq 0,$ the GPSWFs are solutions of an energy maximization problem
 over a generalized Paley-Wiener space and with respect to  certain weighted norms. As  important consequences of this characterization, 
 we get a monotonicity result  the sequence ${\displaystyle \lambda_n^{(\alpha)}(c)=\frac{c}{2\pi} |\mu_n^{(\alpha)}(c)|^2}$  with respect to the parameter $\alpha.$
Moreover, by using the results of \cite{Bonami-Karoui3}, one gets 
a better understanding of the behaviour and the super-exponential decay rate of the $(\lambda_n^{(\alpha)}(c))_{n\geq 0}.$
 In the second part, we show that the GPSWFs are well adapted for the approximation of functions from the classical Paley-Wiener space $B_c$ as well as of almost $c-$band-limited functions.
 
 \subsection{GPSWFs as solutions of an energy maximization problem and consequences.}
 
 We  recall that the starting point of the theory of the classical  PSWFs (corresponding to the GPSWFs with $\alpha=0)$ is the solution of the following energy maximization problem, see \cite{Slepian1}
 \begin{equation}\label{energymax1}
 f = \arg \max_{f \in B_c} \frac{\| f\|_{L^2(I)}^2}{\| f\|_{L^2(\mathbb R)}^2}=\arg \max_{f \in B_c} 2\pi \frac{\| f\|_{L^2(I)}^2}{\|\widehat f\|_{L^2(\mathbb R)}^2},
 \end{equation}
 where $B_c$ is the Paley-Wiener of $c-$band-limited functions given by 
\begin{equation}\label{Bc}
B_c=\{ f\in L^2(\mathbb R),\,\, \mbox{Support } \widehat f\subseteq [-c,c]\}.
 \end{equation} 
More precisely, it has been shown in \cite{Slepian1} that from $B_c,$ $\psi_{0,c}^{(0)} $ is the most 
concentrated function in $I=[-1,1]$ with the largest energy concentration ratio $0<\lambda_0^{(0)}(c)<1.$ Moreover, for any integer $n\geq 1,$
$\psi_{n,c}^{(0)} $ is the most concentrated function from $B_c$ which is orthogonal to the previous $\psi_{i,c}^{(0)},\, 0\leq i\leq n-1.$ The orthogonality is with respect to the two usual inner products of $L^2(I)$ and $L^2(\mathbb R).$ As it will be seen, the extension of the previous 
characterization of the PSWFs to the more general case of the GPSWFs  provides us with a better understanding of the behaviour and the super-exponential decay rate of the eigenvalues $(\lambda_n^{(\alpha)}(c))_{n\geq 0}.$ For $\alpha > 0,$ we define the restricted Paley-Wiener space of weighted $c-$band-limited functions by 
\begin{equation}\label{GBc}
B_c^{(\alpha)}=\{ f\in L^2(\mathbb R),\,\,  \mbox{Support } \widehat f\subseteq [-c,c],\, \, \widehat f\in
L^2\big((-c,c), \omega_{- \alpha}(\frac{\cdot}{c})\big)\}.
 \end{equation} 
Here, $L^2\big((-c,c), \omega_{- \alpha}(\frac{\cdot}{c})\big)$ is the weighted $L^2(-c,c)-$space with norm given by 
$$\| f\|^2_{L^2\big((-c,c), \omega_{- \alpha}(\frac{\cdot}{c})\big)}= \int_{-c}^c |f(t)|^2 \omega_{-\alpha}\left(\frac{t}{c}\right)\, dt.$$
 Note that when $\alpha=0,$ the restricted Paley-Wiener space $B_c^{(0)}$ is reduced to the usual space $B_c.$ Also, since for any $\alpha \geq \alpha',$ $\widehat f\in
L^2\big((-c,c), \omega_{- \alpha'}(\frac{\cdot}{c})\big)$ implies that $\widehat f\in
L^2\big((-c,c), \omega_{- \alpha}(\frac{\cdot}{c})\big)$ then one gets 
\begin{equation}
\label{inclusionPaley}
B_c^{(\alpha)} \subseteq B_c^{(\alpha')},\qquad \forall\, \alpha \geq \alpha' \geq 0.
\end{equation}

\begin{remark} We give an example of a function from a restricted Paley-Wiener space. 
 If $c>1,$ then it has been shown in \cite{Karoui1} that the function $${\displaystyle 
\eta(x)= \frac{\sin(c x-x/2)}{x^4}\left(2\sin\left(\frac{x}{2}\right)-x \cos\left(\frac{x}{2}\right)\right),\,\, x\in \mathbb R}$$ is a $c-$band-limited function.
Moreover, its Fourier transform is an even function given by 
$$\widehat \eta(-\xi)=\widehat \eta(\xi)=\left\{\begin{array}{ll} 1 & \mbox{ if \ \ } \xi\in [0,c-1]\\
f(\xi+2-c) &\mbox{ if \ \ } \xi \in [c-1,c]\\
0 & \mbox{ if \ \ } \xi \geq c.
\end{array}
\right. $$
where
$$f(-x)=f(x)=\left\{\begin{array}{ll} 1 & \mbox{ if \ \ } x \in [0,1]\\
-4+12 x - 9 x^2+2 x^3 &\mbox{ if \ \ } x \in [1,2]\\
0 & \mbox{ if \ \ } x \geq 2.
\end{array}
\right. $$
Since for $\xi\in [c-1,c],$ $\widehat \eta(\xi)= (2\xi-2c+3) (\xi-c)^2,$ then it is easy to see that $\eta$ belongs to the 
restricted Paley-Wiener space  $B_c^{\alpha}$ for any 
$0\leq \alpha <5.$
\end{remark}
 
The generalized maximization problem is formulated as follows. We note that from \eqref{fourier_Jacobi2} with $k=0,$ one gets the finite Fourier transform of the weight function $\omega_{\alpha},$ given by
\begin{equation}
\label{fourier_weight}
\int_{-1}^1 e^{ixy} \omega_{\alpha}(y)\, dy = \sqrt{\pi} 2^{\alpha+1/2}\Gamma(\alpha+1) \frac{J_{\alpha+1/2}(x)}{x^{\alpha+1/2}}= 
\mathcal K_{\alpha}(x),\quad x\in \mathbb R.
\end{equation}
Next, if $f\in B_c^{(\alpha)},$  then $\widehat f (x)= g(x) \omega_{\alpha}\left(\frac{x}{c}\right)$ with some ${\displaystyle g\in 
L^2\big((-c,c), \omega_{ \alpha}(\frac{\cdot}{c})\big).}$ By using the inverse Fourier transform, one gets 
\begin{eqnarray}\label{energymax2}
\frac{\|f\|^2_{L^2(I,\omega_{\alpha})}}{\|\widehat f\|^2_{L^2(\omega_{-\alpha}(\frac{\cdot}{c}))}}&=& 
\frac{1}{\|\widehat f\|^2_{L^2(\omega_{-\alpha}(\frac{\cdot}{c}))}} \int_{-1}^1 f(t)\cdot \overline{f(t)} \omega_{\alpha}(t)\, dt \nonumber\\
&=&\frac{1}{\|\widehat f\|^2_{L^2(\omega_{-\alpha}(\frac{\cdot}{c}))}} \frac{1}{4\pi^2}\int_{-1}^1 \int_{-c}^c e^{it y}\widehat f(y) \, dy\cdot  \int_{-c}^c e^{-it x} \overline{\widehat f(x)} \, dx \omega_{\alpha}(t)\, dt  \nonumber\\
&=& \frac{1}{\|\widehat f\|^2_{L^2(\omega_{-\alpha}(\frac{\cdot}{c}))}} \frac{1}{4\pi^2}\int_{-c}^c \int_{-c}^c
\left(\int_{-1}^1
e^{it(y-x)} \omega_{\alpha}(t)\, dt\right)\widehat f(y) \, dy\cdot  \int_{-c}^c e^{-it x} \overline{\widehat f(x)} \, dx \omega_{\alpha}(t)\, dt \nonumber\\
&=&\frac{1}{\|\widehat f\|^2_{L^2(\omega_{-\alpha}(\frac{\cdot}{c}))}} \frac{1}{4\pi^2} \int_{-c}^c \int_{-c}^c \mathcal K_{\alpha}(y-x) g(y) 
\omega_{\alpha}\left(\frac{y}{c}\right) \, dy \cdot \overline{g(x)} \omega_{\alpha}\left(\frac{x}{c}\right) \, dx \nonumber\\
&=&\frac{1}{4\pi^2} \frac{ {\displaystyle\int_{-c}^c \left[ \int_{-c}^c \mathcal K_{\alpha}(y-x) g(y) 
\omega_{\alpha}\left(\frac{y}{c}\right) \, dy \right] \cdot \overline{g(x)} \omega_{\alpha}\left(\frac{x}{c}\right) \, dx}}{{\displaystyle \int_{-c}^c g(x) \overline{g(x)}
\omega_{\alpha}\left(\frac{x}{c}\right) \, dx}}.
\end{eqnarray}
Here, $\mathcal K_{\alpha}$ is as given by \eqref{fourier_weight}. Note that since the compact integral operator 
$\mathcal Q^{\alpha}$ defined on ${\displaystyle 
L^2(\omega_{\alpha}(\frac{\cdot}{c}))}$ by 
$$ \mathcal Q^{\alpha} g(x) = \frac{1}{4\pi^2} \int_{-c}^c \mathcal K_{\alpha}(y-x) g(y) \omega_{\alpha}\left(\frac{y}{c}\right) \, dy$$
has a symmetric kernel, then it is well known that in this case, 
${\displaystyle \max_{f\in B^{\alpha}_c} 2 \pi \frac{\|f\|^2_{L^2_{\omega_{\alpha}}(I)}}{\|\widehat f\|^2_{L^2(\omega_{-\alpha}(\frac{\cdot}{c}))}} }$
is attained at the eigenfunction of $2\pi \mathcal Q^{\alpha},$ associated with the largest eigenvalue. Hence, by using a trivial change of variable and functions, the generalized energy maximization problem is
reduced to the solution of the following eigenproblem
\begin{equation}
\mathcal Q_c^{(\alpha)} G (x) = \int_{-1}^1 \frac{c}{2 \pi}\mathcal K_{\alpha}(c(x-y)) G(y) \omega_{\alpha}(y) \, dy = \lambda G(x),\quad x\in [-1,1].
\end{equation} 
On the other hand, it has been shown in \cite{Wang2} that the kernel $\mathcal K_{\alpha}(c(x-y))$ is nothing but the kernel of the composition  operators $\mathcal F_c^{{\alpha}*} \circ \mathcal F_c^{\alpha}.$  Hence, the operators $\mathcal Q_c^{\alpha}$ and $\mathcal F_c^{\alpha}$ have the same  eigenfunctions, given by the GPSWFs, $\psi_{n,c}^{(\alpha)},$ and associated to the respective eigenvalues $\lambda_n^{(\alpha )},$ $\mu_n^{(\alpha)}(c).$
These eigenvalues are  related
to each others by the following rule
\begin{equation}
\lambda_n^{(\alpha)}(c)= \frac{c}{2\pi } |\mu_n^{(\alpha)}(c)|^2,\,\, n\geq 0.
\end{equation}
 Since from \eqref{inclusionPaley}, if $0\leq \alpha'\leq \alpha,$ then we have $B_c^{(\alpha)} \subseteq B_c^{(\alpha')}$ and since for $f\in B_c^{(\alpha)},$ then we have 
$$\|f\|^2_{L^2(I,\omega_{\alpha})}\leq \|f\|^2_{L^2(I,\omega_{\alpha'})},\quad 
\|\widehat f\|^2_{L^2(\omega_{- \alpha}(\frac{\cdot}{c}))} \geq \|\widehat f\|^2_{L^2(\omega_{- \alpha'}(\frac{\cdot}{c}))}.$$
Hence, we have
$$\lambda_0^{(\alpha)}(c)=\sup_{f\in B_c^{(\alpha)}} \frac{\|f\|^2_{L^2(I,\omega_{\alpha})}}
{\|\widehat f\|^2_{L^2(\omega_{- \alpha}(\frac{\cdot}{c}))}} \leq \sup_{f\in B_c^{(\alpha)}} \frac{\|f\|^2_{L^2(I,\omega_{\alpha'})}}{\|\widehat f\|^2_{L^2(\omega_{- \alpha'}(\frac{\cdot}{c}))}}\leq \sup_{f\in B_c^{(\alpha')}} \frac{\|f\|^2_{L^2(I,\omega_{\alpha'})}}{\|\widehat f\|^2_{L^2(\omega_{- \alpha'}(\frac{\cdot}{c}))}} =\lambda_0^{(\alpha')}(c).$$
More generally, for an integer $n\geq 1,$ let $\psi_{0,c}^{(\alpha')},\ldots,\psi_{n-1,c}^{(\alpha')}$ be the first most concentrated
GPSWFs, associated with the respective eigenvalues $\lambda_0^{(\alpha')}(c)>\lambda_1^{(\alpha')}(c)>\cdots>\lambda_{n-1}^{(\alpha')}(c).$
 Note that the previous strict inequalities are due to the fact that these eigenvalues are simple, see \cite{Wang2}. By combining the previous formulation 
 of the energy maximization problem and  the well known Min-Max principle for eigenvalues of compact operator, one concludes that if $S_n,\, H_n$ stand for an arbitrary subspace of dimension $n$ of $B_c^{(\alpha)},$
 and $B_c^{(\alpha')},$ respectively, then we have 
 \begin{eqnarray*}
 \lambda_n^{(\alpha)}(c)&=&\max_{S_n\subset B_c^{(\alpha)}} \min_{\psi \in S_n} 
 \frac{\|\psi \|^2_{L^2(I,\omega_{\alpha})}}{\|\widehat \psi \|^2_{L^2(\omega_{- \alpha}(\frac{\cdot}{c}))}}
 \leq \max_{S_n\subset B_c^{(\alpha')}} \min_{\psi \in S_n} 
  \frac{\|\psi \|^2_{L^2(I, \omega_{\alpha'}))}}{\|\widehat \psi \|^2_{L^2(\omega_{- \alpha'}(\frac{\cdot}{c}))}}\\
 &\leq& \max_{H_n\subset B_c^{(\alpha')}} \min_{\psi \in H_n} 
  \frac{\|\psi \|^2_{L^2(I,\omega_{\alpha'})}}{\|\widehat \psi \|^2_{L^2(\omega_{- \alpha'}(\frac{\cdot}{c}))}} = \lambda_n^{(\alpha')}(c).
 \end{eqnarray*}

We have just proved the following theorem giving the monotony of the eigenvalues $\lambda_n^{(\alpha)}(c)$  with respect to the parameter $\alpha.$ 

\begin{theorem}\label{monotony_alpha}
For a given real number  $c>0,$  and an integer $n\geq 0,$  we have 
\begin{equation}
\lambda_n^{(\alpha)}(c) \leq \lambda_n^{(\alpha')}(c),\quad \forall\, \alpha \geq \alpha' \geq 0. 
\end{equation}
\end{theorem}

It is important to mention that a super-exponential 
  decay rate of the sequence  $(\lambda_n^{(\alpha)}(c))_n$ as well as an estimate of the location of the plunge region, where  the fast decay starts
  are important consequences of the previous proposition. These two results
 follow directly from the results given in \cite{Bonami-Karoui3}, where an
  explicit formula for estimating the $\lambda_n^{(0)}(c)$ has been developed. This explicit formula enjoys with a surprising accuracy as soon as $n$ reaches or goes beyond the plunge region around the value $n_c = \frac{2c}{\pi}.$ Also, it proves that the exact asymptotic  super-exponential decay rate is given by 
 the quantity ${\displaystyle e^{-2n \log\left(\frac{4 n}{ec}\right)}.}$ From the previous theorem  with $\alpha'=0,$ one concludes that
 for any $\alpha >0,$  the sequence  $(\lambda_n^{(\alpha)}(c))_n$ has a super-exponential decay rate, bounded above by the decay rate of $(\lambda_n^{(0)}(c))_n.$ Moreover, the fast decay of these $(\lambda_n^{(\alpha)}(c))_n$ starts around  
 $n_c = \frac{2c}{\pi}.$ In  the numerical results section, we give different tests that illustrate these precise  behaviours of $(\lambda_n^{(\alpha)}(c))_n.$

 \subsection{Approximation of Band-limited functions by the GPSWFs.}
In this paragraph, we first show that when restricted to the interval $I,$ the GPSWFs $\psi_{n,c}^{(\alpha)}$ are  well adapted for the approximation  of functions
from the usual Paley-Wiener space $B_c.$ As a result, we check that the GPSWFs are also well adapted for the approximation of almost band-limited functions. This  type of functions have been defined in \cite{Landau} as follows. 

\begin{definition}
 Let $\Omega=[-c, c],$ then a 
function  $f$ is said to be $\epsilon_{\Omega}-$band-limited in $\Omega$ if
$$\frac 1{2\pi}\int_{|\xi|> c} |\widehat f(\xi)|^2\, d\xi \leq
\epsilon_{\Omega}^2.$$
\end{definition}

 \begin{proposition}\label{band-limited} Let  $c>0,\,\, \alpha \geq 0$ be two real numbers and let $f\in B_c.$
 For any positive integer $N> \frac{2c}{\pi},$ let 
 $$S_N(f)(x) = \sum_{k=0}^{N} <f,\psi_{k,c}^{(\alpha)}>_{L^2(I,\omega_{\alpha})} \psi_{k,c}^{(\alpha)}(x).$$
 Then, we have
\begin{equation}\label{approx1}
 \left(\int_{-1}^{1}|f(t)- S_N f(t)|^2 \omega_{\alpha}(t) dt\right)^{1/2}\leq   C_1 \sqrt{\lambda_N^{(\alpha)}(c)}\, (\chi_N(c))^{(1+\alpha)/2}\, \| f\|_{L^2(\mathbb R)},
\end{equation}
and
\begin{equation}\label{approx2}
\sup_{x\in [-1,1]}|f(x)- S_N f(x)|\leq   C_1 \sqrt{\lambda_N^{(\alpha)}(c)}\, (\chi_N(c))^{1+\alpha/2}\, \| f\|_{L^2(\mathbb R)}.
\end{equation}
for some uniform constant $C_1$ depending only on $\alpha.$
\end{proposition}

\noindent
{\bf Proof:} We first note that since $\mathcal B=\{\psi_{n,c}^{(\alpha)},\,\, n\geq 0\}$ is an orthonormal basis of 
$L^2(I,\omega_{\alpha}),$ and since $\chi_I f\in L^2(I),$ where $\chi_I$ denotes the characteristic function, then we have
\begin{equation}
\label{expansion3}
f(x)= \sum_{k\geq 0} <f,\psi_{k,c}^{(\alpha)}>_{L^2(I,\omega_{\alpha})} \psi_{k,c}^{(\alpha)}(x),\quad a.e.\,\,  x\in I.
\end{equation}
On the other hand, since $f\in B_c,$ then $f\in \mathcal C^{\infty}(\mathbb R)\cap L^2(\mathbb R).$ In particular, from the inverse Fourier transform,
and by using the fact that $f\in B_c,$ we have
\begin{equation}\label{Eqq1}
f(x)= \frac{1}{2\pi} \int_{\mathbb R} e^{i x y}\widehat f(y) dy=\frac{1}{2\pi} \int_{-c}^c e^{i x y}\widehat f(y) dy=\frac{c}{2\pi }\int_{-1}^1 e^{ictx} \widehat f(ct) \, dt,\quad \forall\, x\in [-1,1].
\end{equation}
Consequently, for any integer $k\geq 0,$ we have
\begin{eqnarray}
\lefteqn{|<f,\psi_{k,c}^{(\alpha)}>_{L^2(I,\omega_{\alpha})}|=\left|\int_{-1}^1 f(x) \psi_{k,c}^{(\alpha)}(x) \omega_{\alpha}(x)\, dx\right|}\\
& =&
\frac{c}{2\pi }\left|\int_{-1}^1\left(\int_{-1}^1 e^{ictx} \widehat f(ct) \, dt\right) \psi_{n,c}^{(\alpha)}(x) \omega_{\alpha}(x)\, dx\right|= \frac{c}{2\pi} \left|\int_{-1}^1 \widehat f(ct) \left(\int_{-1}^1 e^{ictx} \psi_{k,c}^{(\alpha)}(x) \omega_{\alpha}(x)\, dx\right)\, dt\right|\nonumber\\
&=& \frac{c}{2\pi}\, |\mu_k^{(\alpha)}(c)| \left|\int_{-1}^1 \widehat f(ct) \psi_{n,c}^{(\alpha)}(t)\, dt\right| \leq  |\mu_k^{(\alpha)}(c)| \frac{c}{2\pi} \sup_{t\in [-1,1]}|\psi_{n,c}^{(\alpha)}(t)| \left(\int_{-1}^1 |\widehat f(ct)|^2\, dt\right)^{1/2}\nonumber \\
&\leq & \frac{C_{\alpha}}{\sqrt{2}} \sqrt{\frac{c}{2\pi}} |\mu_k^{(\alpha)}(c)| (\chi_k(c))^{(1+\alpha)/2} \|f\|_{L^2(\mathbb R)}= \frac{C_{\alpha}}{\sqrt{2}} \sqrt{\lambda_N^{(\alpha)}(c)}\, (\chi_N(c))^{(1+\alpha)/2}\, \| f\|_{L^2(\mathbb R)}.   
\end{eqnarray}
Here, $C_{\alpha}$ is as given by \eqref{maxpsi1}. 
The last inequality follows from Plancherel formula and the bound over $I$ of $|\psi_{n,c}^{(\alpha)}(t)|,$ we have given in \eqref{maxpsi2}.
On the other hand, by using the previous inequality,  together with the super-exponential decay rate of the $|\mu_n^{(\alpha)}(c)|,$ given by \eqref{decaymu_n}, 
as well as the Parseval's  equality $${\displaystyle \|f-S_Nf\|^2_{L^2(I,\omega_{\alpha})}=\sum_{k\geq N+1}|<f,\psi_{k,c}^{(\alpha)}>_{L^2(I,\omega_{\alpha})}|^2},$$ one can easily get \eqref{approx1}.
Finally, to get \eqref{approx2}, it suffices to combine the previous inequality, \eqref{decaymu_n} as well as the upper bound of $|\psi_{n,c}^{(\alpha)}(t)|.$

As a consequence of the previous result, we have the following corollary concerning the quality of approximation
of almost band-limited functions by the GPSWFs.

\begin{corollary}\label{approx_almostband}
Let $f\in L^2(\mathbb R)$ be an 
$\epsilon_{\Omega}-$band-limited in $\Omega=[-c, +c]$ and let $\alpha\geq 0,$  then
for any positive integer $N \geq \frac{2c}{\pi},$ we have
\begin{equation}\label{eqqq6.0}
\left(\int_{-1}^{1}|f(t)- S_N f(t)|^2 \omega_{\alpha}(t) dt \right)^{1/2}\leq  \epsilon_{\Omega}+C_1 \sqrt{\lambda_N^{(\alpha)}(c)}\, (\chi_N(c))^{(1+\alpha)/2}\, \| f\|_{L^2(\mathbb R)}
 \end{equation}
 where the constant $C_1$ depends only on $\alpha.$
\end{corollary}

\noindent
{\bf Proof:} It suffices to consider the band-limiting operator $\pi_{\Omega}$ defined by:
$$\pi_{\Omega}(f)(x)=\frac{1}{2\pi}\int_{\Omega} e^{i x
\omega}\widehat{f}(\omega)\, d\omega.$$
Since $\pi_{\Omega} f \in B_c,$ $\| (f-\pi_{\Omega} f)- S_N(f-\pi_{\Omega} f)\|_{L^2(I)}\leq \|f-\pi_{\Omega}\|_{L^2(\mathbb R)}\leq \epsilon_{\Omega}$ and $\| \pi_{\Omega} f \|_{L^2(\mathbb R)}\leq \| f \|_{L^2(\mathbb R)},$ then by applying the result of the previous proposition to $\pi_{\Omega} f,$ one gets 
\begin{eqnarray*}
\|f- S_N(f)\|_{L^2(I,\omega_{\alpha})} &\leq & \|(f-\pi_{\Omega} f)- S_N (f-\pi_{\Omega} f)\|_{L^2(I,\omega_{\alpha})} + \|\pi_{\Omega} f- S_N (\pi_{\Omega} f)\|_{L^2_{\omega_{\alpha}(I)}} \\
&\leq & \epsilon_{\Omega}+ C_1 \sqrt{\lambda_N^{(\alpha)}(c)}\, (\chi_N(c))^{(1+\alpha)/2}\, \|\pi_{\Omega} f\|_{L^2(\mathbb R)} \\
&\leq &
 \epsilon_{\Omega}+ C_1 \sqrt{\lambda_N^{(\alpha)}(c)}\, (\chi_N(c))^{(1+\alpha)/2} \, \| f\|_{L^2(\mathbb R)}.
\end{eqnarray*}

\section{Numerical results.}

In this section, we give three examples that illustrate the different results of this work. The first example deals with the computation and the analytic extension of the GPSWFs.\\

\noindent
{\bf Example 1:} In this  example, we   give different  numerical tests that illustrate the construction scheme of the GPSWFs $\psi_{n,c}^{(\alpha)}.$  For this purpose, we have considered the values  $\alpha=0.5$ and  $c=5\pi.$ Then, we have computed the different Jacobi expansion coefficients via the scheme of section 2, by solving the eigensystem \eqref{eigensystem1}, truncated to the order $N=90.$ Figures 1(a) show the graphs of the  $\psi_{n,c}^{(\alpha)}$ for the different values of $n=0, 5, 15.$ Note that these graphs illustrate some of the provided  properties of the $\psi_{n,c}^{(\alpha)}.$ Also, we have used formula \eqref{expansion3} and computed the analytic extensions of the previous GPSWFs. The graphs of these extensions are given by Figure 1(b) and 1(c). Note that as predicted by the
characterisation of the GPSWFs as solutions of the energy maximization problem, for the  values of 
$n\leq 2c\pi,$ the  $\psi_{n,c}^{(\alpha)}$ are concentrated on $I,$ whereas for $n> 2c/\pi,$ they are concentrated on $\mathbb R\setminus I.$ \\

 \begin{figure}[h]\hspace*{0.05cm}
 {\includegraphics[width=15.05cm,height=5.5cm]{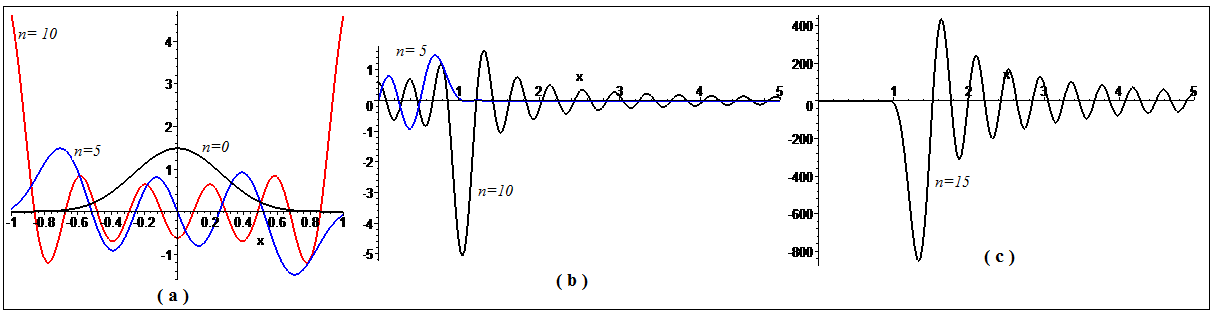}}
 \vskip -0.5cm\hspace*{1cm} \caption{{\large (a) Graphs of $\psi_{n,c}^{\alpha}$ with $c=5\pi$ and $\alpha=0.5,$\,\,\, 
 (b) and (c) graphs of some analytic extensions of the  $\psi_{n,c}^{\alpha}.$ } }
 \end{figure}

\noindent
{\bf Example 2:} In this example, we illustrate the important result given by Theorem \ref{monotony_alpha}, concerning the 
monotonicity with respect to the parameter $\alpha$ of the sequence ${\displaystyle \lambda_n^{(\alpha)}(c)=\frac{c}{2\pi} |\mu_n^{(\alpha)}(c)|^2.}$
For this purpose, we have  used formula \eqref{mu_n} and computed highly accurate values of $\mu_n^{(\alpha)}(c)$ and consequently of 
$\lambda_n^{(\alpha)}(c)$ with 
$c=10 \pi$ and with different values of $\alpha=0,\, 0.5,\, 1.5.$   Note that as predicted 
by Theorem \ref{monotony_alpha}, the sequence $\lambda_n^{(\alpha)}(c)$ is decreasing with respect to $\alpha.$  
The graphs of the $\lambda_n^{(\alpha)}(c)$ as well as $\log(\lambda_n^{(\alpha)}(c)) $ are given by Figure 2(a),and 2(b), respectively.\\

\begin{figure}[h]\hspace*{0.05cm}
{\includegraphics[width=15.05cm,height=5.5cm]{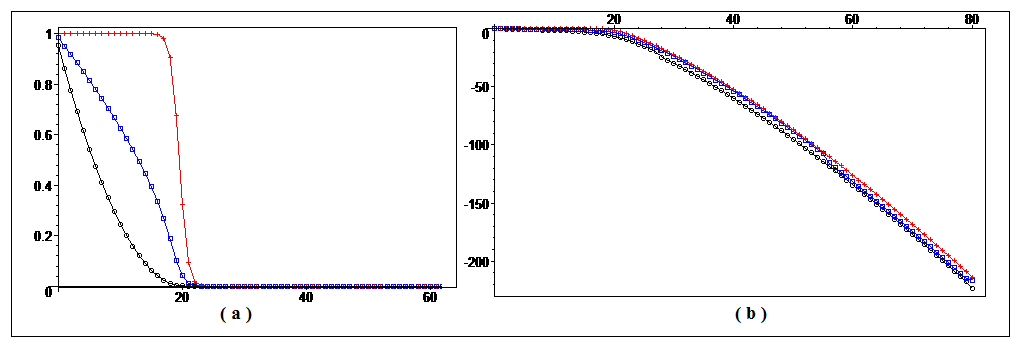}}
\vskip -0.5cm\hspace*{1cm} \caption{{\large (a) Graphs of the $\lambda_n^{(\alpha)}(c)$ for $c=10\pi,$ $\alpha=1.5$ (circles),
$\alpha=0.5$ (boxes) and $\alpha=0$ (crosses),\,\, (b) graphs of  $\log(\lambda_n^{(\alpha)}(c)).$}}
\end{figure}

\noindent
{\bf Example 3:} In this last example, we illustrate the quality of approximation over $I$ of band-limited and almost
 band-limited functions, by the GPSWFs. For this purpose, we have first considered the value of $\alpha=0.5$ and the $c-$band-limited function ${\displaystyle f(x)=\frac{\sin(cx)}{cx}}$ with $c=50.$ By computing 
 the projections $S_N(f),$ with $N=32$ and $N=40,$ we found that 
 $$\sup_{x\in [-1,1]}|f(x)-S_{32}(f)(x)|\approx 2.22\, 10^{-2},\quad  \sup_{x\in [-1,1]}|f(x)-S_{40}(f)(x)|\approx 4.80\, 10^{-6}.$$
As predicted by proposition \ref{band-limited},  the drastic improvement in the previous approximation errors
is due to the fact that the second value of $N=40$ lies after  the plunge region of the eigenvalues
$\lambda_n^{(\alpha)}(c),$ which is not the case for the first value of $N=32.$ \\ 
 
 Next, to illustrate the approximation of almost band-limited functions by the GPSWFs, we have considered the
 Weierstrass function
 \begin{equation}\label{eq5.1}
 W_s(x)= \sum_{k\geq 0} \frac{\cos(2^{k}x)}{2^{ks}},\quad -1\leq
 x\leq 1.
 \end{equation}
 It is well known that  $W_s \in H^{s-\epsilon}(I),\,\forall \epsilon <  s,\, s  >0.$
 One may consider $W_s$ as a restriction over $I$ of a function $W\in H^{s-\epsilon}(\mathbb R).$
Note that  if $f\in H^s(\R)$ with $s>0$,  then
 $$
 \int_{|\xi|> c}|\widehat{f}(\xi)|^2\,\mbox{d}\xi \leq 
 \int_{|\xi|> c}\frac{(1+|\xi|)^{2s}}{(1+|\xi|)^{2s}}|\widehat{f}(\xi)|^2\,\mbox{d}\xi
 \leq \frac{1}{(1+c)^{2s}} \norm{f}_{H^s(\mathbb R)}^2.
 $$
 That is  $f$ is $\frac{1}{(1+c)^{s}}\norm{f}_{H^s}-$almost band-limited to
 $[-c,c]$.
  Note that in \cite{Bonami-Karoui4}, we have 
 used the previous function to illustrate the quality of approximation by the classical PSWFs. In this example, we push forward this quality of 
 approximation to the GPSWFs. For this purpose, we have considered the value of $\alpha=0.5$ and the  two couples of $(c,N)=(50,60), (100, 90).$   Then we have computed the 
 associated projection $S_N( W_s),$ for the  value of $s=1.$ Note  that thanks to \eqref{fourier_Jacobi2}, the different expansion coefficients 
 ${\displaystyle C_n(W_s)=\int_{-1}^1 W_s(y) \psi_{n,c}^{(\alpha)}(y) \omega_{\alpha}(y)\, dy}$ are computed exactly.
 In fact since $W_s$ is an even function, and from the Jacobi series expansion of $\psi_{n,c}^{(\alpha)},$ the computation of the  $C_n(W_s)$ is restricted to the even indexed coefficients and consequently to the computation of the different inner products with Jacobi polynomials of even degrees. More precisely, we have
 \begin{eqnarray*}
 C_{2 m}(W_s) &=& \sum_{l\geq 0} \beta_{2 l}^{2m} \sum_{k\geq 0}\frac{1}{2^{ks}}\int_{-1}^1 \cos(2^{k}x) \widetilde P_{2l}^{(\alpha,\alpha)}(y)\, \omega_{\alpha}(y)\, dy\\
 &=& \sqrt{\pi} 2^{\alpha+1/2} \sum_{l\geq 0} (-1)^l \beta_{2 l}^{2m} \frac{\Gamma(2l+\alpha+1)}{\sqrt{h_{2l}} (2l)!} \sum_{k\geq 0}  2^{-k(s+\alpha+1/2)} J_{2l+\alpha+1/2}(2^k).
  \end{eqnarray*}
  \begin{figure}[h]\hspace*{0.05cm}
  {\includegraphics[width=17.05cm,height=5.5cm]{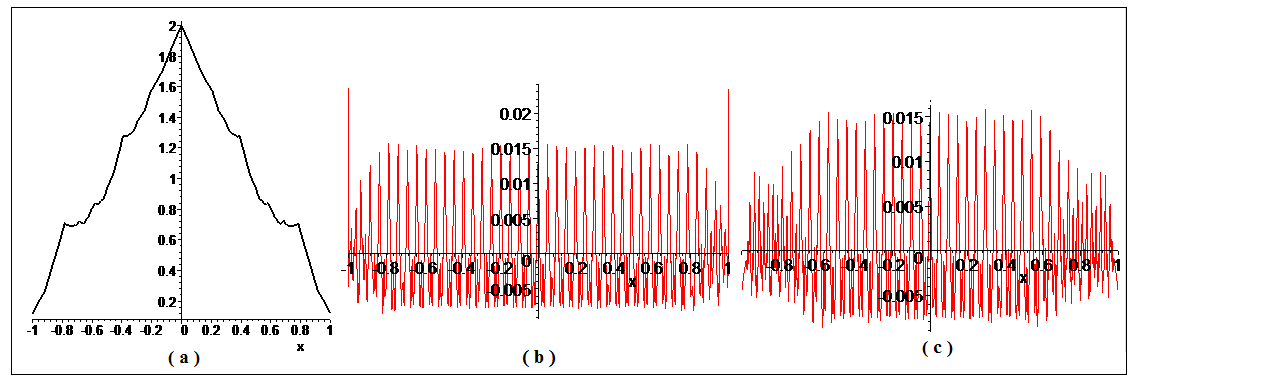}}
  \vskip -0.5cm\hspace*{1cm}  \caption{{\large
 (a) Graph of the Weierstarss function $W_1(x),$ 
   (b) Graph of the approximation error $W_1(x)- S_N(W_1)(x)$ with $\alpha=0.5,$ $(c,N)=(50,60),$  (c) same as (b) with $(c,N)=(100,90).$}}
     \end{figure}
The graph of $W_1$ is given by Figure 3(a), whereas the graphs of the  approximation errors $W_1(x)-S_N(W_1)(x)$ corresponding to the two couples
$(c,N)=(50,60), (100, 90)$ are given by Figure 3(b) and 3(c), respectively. Note that as predicted by the theoretical results of section 4, the 
approximation error decreases as $c^{-s},$ whenever the truncation order $N$ lies beyond the plunge of the 
$(\lambda_n^{(\alpha)}(c))_n.$ In this case, the extra error factor given by $\sqrt{\lambda_N^{(\alpha)}(c)} \chi_N^{1/2+\alpha}$ can be neglected comparing the factor $c^{-s}.$\\

\noindent
{\bf Acknowledgement:} The authors thank very much the anonymous referees for their valuable comments that
have improved this work. 

\vskip 0.5cm
\noindent
{\large {\bf Appendix:}} Proof of Theorem 1.\\

\noindent
The proof is divided into three steps. To alleviate notations of this proof, we will simply denote 
$\psi_{n,c}^{(\alpha)}$ and $\chi_n(c)$ by $\psi_{n,c}$ and $\chi_n,$ respectively.\\

\noindent 
{\it First step:} We prove that for for any positive integer $j$ with  $j(j+2\alpha+1)\leq \chi_n$, all moments $\int_{-1}^1 y^j \psi_{n,c}(y)\, dy$ are non negative and 
\begin{equation}
\label{momentspsi}
0\leq\int_{-1}^1 y^{j} \psi_{n,c}(y)\, \omega_{\alpha}(y)\, dy \leq \sqrt{1+\alpha}\left(\frac{1}{q}\right)^{j} |\mu_n^{(\alpha)}(c)|.
\end{equation}
To this end, we first check that for any  integer $ k\geq 0$ satisfying $k(k+2\alpha+1)\leq \chi_n,$ we have
\begin{equation}\label{ineqderivatives}
\left|\psi_{n,c}^{(k)}(0)\right|\leq  (\sqrt{\chi_n})^k \sqrt{1+\alpha}.
\end{equation}
It suffices to prove that $ m_k=\frac{|\psi_{n,c}^{(k)}(0)|}{\sqrt{\chi_n}^k} \leq \sqrt{1+\alpha}.$  
From the parity of $\psi_{n,c},$ we need  only to consider  derivatives of even or odd order. We 
assume that $n=2l$ is even. The case where $n$ is odd is done in a similar manner.
Note that  for a fixed $ n,$ $ \psi_{n,c}^{(2l)}(0) $ has alternating signs, that is  $ \psi_{n,c}^{(k)}(0)  \psi_{n,c}^{(k-2)}(0)<0 $ 
In fact, for $k=0,$ we have $ \psi_{n,c}(0)\psi_{n,c}^{(2)}(0)=-\chi_n\psi_{n,c}(0)^2<0.$ By induction,  we assume that 
$ \psi_{n,c}^{(k)}(0) \psi_{n,c}^{(k-2)}<0 .$ As it is done in \cite{Bonami-Karoui3}, 
 we have
\begin{equation}
\psi_{n,c}^{(k+2)}(0) \psi_{n,c}^{(k)}(0)=\Big( k(k+1+2\alpha)-\chi_n \Big) \psi_{n,c}^{(k)}(0)^2+k(k-1)c^2 \psi_{n,c}^{(k-2)}(0)\psi_{n,c}^{(k)}(0).
\end{equation}
By using the induction hypothesis as well as the fact that $ k(k+1+\alpha+\beta) \leq \chi_n,$ one concludes that the induction assumption holds for the order $k.$  Consequently, we have
\begin{equation}
|\psi_{n,c}^{(k+2)}(0)|=\Big( \chi_n-k(k+1+2\alpha)\Big) |\psi_{n,c}^{(k)}(0)|+k(k-1)c^2| \psi_{n,c}^{(k-2)}(0)|.
\end{equation}
The previous equality implies that 
\begin{equation}
m_{k+2}=\Big(1-\frac{k(k+1+2\alpha)}{\chi_n}\Big)m_k +k(k-1)\frac{q}{\chi_n}m_{k-2}.
\end{equation}
Hence, for any positive and even integer $k$ with $k(k+2\alpha+1)\leq \chi_n,$ we have 
$m_k\leq m_0\leq \sqrt{1+\alpha}.$ This last inequality follows from \eqref{localest2} with $t=0.$ This proves
the inequality \eqref{ineqderivatives}. Moreover, by taking the $j-$th derivative at zero on both sides of ${\displaystyle \int_{-1}^1 e^{icxy} \psi_{n,c}(y)\omega_{\alpha}(y) dy =\mu_n^{(\alpha)}(c) \psi_{n,c}(x),}$ one gets
\begin{equation}
\label{moments2psi}
\int_{-1}^1 y^j \psi_{n,c}(y)\,\omega_{\alpha}(y) dy= (-i)^j c^{-j} \mu_n^{(\alpha)}(c)\psi_{n}^{(j)}(0).
\end{equation}
Since $\psi_{n}^{(j)}(0)$ and  $\psi_{n}^{(j+2)}(0)$ have opposite signs, then the previous equation implies
that  all moments with even order $j$ with $j(j+2\alpha+1)\leq \chi_n$ have the same sign. The  inequality
 (\ref{momentspsi}) follows from \eqref{ineqderivatives}.\\
 
\noindent
{\it Second Step:}  We show that for all positive  integers $k, n$ with  $k(k+2\alpha+1)+C_{\alpha}\, c^2\leq \chi_n(c)$, we have $\beta_k^n\geq 0.$
Here $C_{\alpha}$ is as given by \eqref{C_alpha}.
The positivity of $\beta_0^n$ (when $n$ is even)  and $\beta_1^n$ (when $n$ is odd) follow from the fact that
\begin{equation}
\label{beta0}
 \beta_0^n=\sqrt{\frac{\Gamma(\alpha+3/2)}{\sqrt{\pi}\Gamma(\alpha+1)}} |\mu_n^{(\alpha)}(c)| |\psi_{n,c}(0)|,\qquad  \beta_1^n=\sqrt{\frac{\Gamma(\alpha+3/2)}{\sqrt{\pi}\Gamma(\alpha+1)}}\sqrt{2\alpha+3}  |\mu_n^{(\alpha)}(c)| |\frac{\psi_n'(0)}{c}|.
\end{equation}
Since the $\beta_k^n$ are given by  (\ref{eigensystem}), then by using the hypothesis of the theorem, we have 
$$ \beta_2^{n} =\frac{2\alpha+3}{c}\sqrt{\frac{2(2\alpha+5)}{2\alpha+2}}\left(\chi_n-\frac{c^2}{2\alpha+3}\right)\beta_0^{n} \geq \beta_0^n,\quad
 \beta_3^{n} =\frac{2\alpha+5}{2c^2}\sqrt{\frac{3(\alpha+1)}{2\alpha+7}}\Big( (2\alpha+2)+\frac{3c^2}{2\alpha+5}\Big)\beta_1^{n}\geq 0.$$
For $j\geq 2$ and by rearranging the system  (\ref{eigensystem}) and using the induction hypothesis $\beta_j^n \geq \beta_{j-2}^n \geq 0,$
one gets 
\begin{equation}\label{ineqq1}
M_{\alpha} c^2 ( \beta_{j+2}^n +\beta_{j-2}^n)\geq
 ( \chi_n(c)-j(j+2\alpha+1)-N_{\alpha} c^2) \beta_j^n, 
\end{equation}
where $M_{\alpha}$ and $N_{\alpha}$ are as given by \eqref{C_alpha}.
If we suppose that $\beta_{j+2} \leq \beta_j^n,$ then from (\ref{ineqq1}), one gets  
\begin{equation}\label{ineq}
2M_{\alpha} c^2  \beta_{j}^n \geq
 ( \chi_n(c)-j(j+2\alpha+1)-N_{\alpha} c^2) \beta_j^n
\end{equation}
which contradicts   the choice of $ C_{\alpha} $ and the fact that $k(k+2\alpha+1)+C_{\alpha}\, c^2\leq \chi_n(c).$
Hence, the induction hypothesis holds for  $\beta_{j+2}^n.$\\

\noindent
{\it Third Step:} We prove \eqref{Decay2beta}. The first inequality follows from  \eqref{beta0} and \eqref{ineqderivatives}.   To prove the second inequality, we recall  that the moments $M_{j,k}$ of the normalized Jacobi polynomials $\widetilde P_k^{(\alpha,\alpha)}$ are given by 
\eqref{moment_Jacobi2} and they are non-negative. Moreover, since ${\displaystyle x^j=\sum_{k=0}^j M_{jk} \widetilde P_k^{(\alpha,\alpha)}(x),}$ then the moments of the $\psi_{n,c}$ are related to the
GPSWFs series expansion coefficients by the following relation
$$\int_{-1}^1 x^j \psi_{n,c}(x)\,\omega_{\alpha}(x)\,  dx = \sum_{k=0}^j M_{j,k}  \beta_k^n.$$
Since from the previous step, we have $\beta_k^n \geq 0,$ for any $0\leq  k\leq j$ and since the $a_{jk}$ are non negative, then the previous equality
implies that
\begin{equation}\label{decay3beta}
\beta_j^n \leq \frac{1}{M_{j,j}} \int_{-1}^1 x^j \psi_{n,c}(x)\,\omega_{\alpha}(x)\,  dx \leq  
\frac{1}{M_{j,j}} \sqrt{1+\alpha}\left(\frac{1}{q}\right)^{j}|\mu_n^{(\alpha)}(c)|.
\end{equation}
The last inequality follows from the result of the first step. Moreover, by using the explicit expression of $M_{j,j},$ given by 
 \eqref{moment_Jacobi2}, together with \eqref{Ineq2},  \eqref{Ineq7}, the decay of the function $\varphi,$ given by 
 \eqref{function}, as well as some straightforward computations, one obtains
\begin{equation}\label{decay4beta}
\frac{1}{M_{j,j}}\leq 2^j \frac{2^{\alpha}}{e^{2\alpha +3/2}}\frac{(3/2)^{3/4}(3/2+2\alpha)^{3/4+\alpha}}{\sqrt{3/2+\alpha}}. 
\end{equation}
Finally, by combining (\ref{decay3beta}) and (\ref{decay4beta}), one gets the second inequality of (\ref{Decay2beta}).$\qquad \Box$

\end{document}